\documentclass[11pt,a4paper]{article}
\usepackage{amssymb}
\usepackage{amsmath}
\usepackage{amsthm}
\usepackage{graphicx}
\usepackage{color}
\usepackage{texdraw}
\usepackage{epsfig}
 \allowdisplaybreaks \numberwithin{equation}{section}
\newtheorem{theorem}{Theorem}[section]
\newtheorem{proposition}[theorem]{Proposition}
\newtheorem{lemma}[theorem]{Lemma}

\newtheorem{example}[theorem]{Example}
\newtheorem{definition}[theorem]{Definition}
\newtheorem{remark}[theorem]{Remark}

\def\neweq#1{\begin{equation}\label{#1}}
\def\endeq{\end{equation}}

 \font \ninemaius=cmcsc10 scaled 900
\font\nineit=cmti9  \font\svfilt=msbm7
\font\ninebf=cmbx9

\font\nineit=cmti9
\font\ninerm=cmr9

\def \div {\mathop {\rm div}\nolimits}
\def \det {\mathop {\rm det}\nolimits}
\def \e {\varepsilon}
\def \osigma {\overline{\sigma}}
\def \u {\overline{u}}
\def \w {\overline \theta _V}
\def \ov {\overline}

\def\qed{{}\hfill{$\square$}\par\medskip}

\def \re {\mathbb R}
\def \svre {\hbox{\svfilt R}}

\def \wt {\widetilde}
\def\Upsilon{\wt{\cal I}}

\begin{document}

\title{A variational method for second order shape derivatives}

\author{Guy BOUCHITT\'E$^\sharp$, \  Ilaria FRAGAL\`A$^*$, \ Ilaria LUCARDESI$^{\sharp\flat}$
\\
{\small $\sharp$ Laboratoire IMATH, Universit\'e de Toulon,
+83957 La Garde Cedex (France)}
\\
{\small $*$ Dipartimento di Matematica, Politecnico di Milano,
Piazza L. da Vinci, 20133 Milano (Italy)}
\\
{\small $\flat$ Present affiliation: Dipartimento di Matematica, SISSA, Trieste (Italy)}
}

\date{}
\maketitle

\begin{abstract} We consider shape  
functionals obtained as minima on Sobolev spaces of classical integrals having smooth and convex densities, 
under mixed Dirichlet-Neumann boundary conditions. 
We propose a new approach for the computation of the second order shape derivative
of such functionals, yielding a general existence and representation theorem. 
In particular, we consider the $p$-torsional rigidity functional for $p \geq 2$.

\end{abstract}



\section{Introduction}
Aim of this paper is to provide a new method for the computation of second order shape derivatives,
which applies to a broad class of shape functionals associated with classical problems in the Calculus of Variations. 
Let us recall that, if $J(\cdot)$ is a functional depending on a subset $\Omega$ of $\re ^n$
and  $V$ is a deformation field in $C^1 (\re^n;\re^n)$, 
the first and second order shape derivatives of $J$ at $\Omega$ in direction $V$, if they exist, are given 
respectively by the limits
$$
J'(\Omega, V):= \displaystyle{\lim_{\e\to 0 } \, \frac{J(\Omega_\e)-J(\Omega)}{\e} }
\, , \qquad 
J''(\Omega,V):=\displaystyle{\lim_{\e\to 0 } 2\, \frac{J(\Omega_\e)-J(\Omega)- \e J'(\Omega, V)}{\e^2} }\,,
$$
where $\Omega _\e$ are the perturbed domains
$
\Omega _\e:= \Psi _\e (\Omega)$,  with $\Psi_\e(x):=x+\e V(x)$.

The literature on  shape derivatives is very wide, and in recent years it has seen a rapid flourishing, also  stimulated by the advances in numerical methods for the determination of optimal shapes. 
We limit ourselves to quote the monographs \cite{DZ, HP, SZ},  
where a lot of bibliographical references in this field can be found.
Concerning in particular second order shape derivatives, their computation is usually irksome, 
but often deserves some efforts because the study of their sign allows to detect whether a critical shape, namely a domain with vanishing first order shape derivative, is actually an optimal one.  As contributions in this direction,  let us mention
without any attempt of completeness the papers \cite{AFM, Bu, DaPi, Ep, Fu, GM, LaPi, NP}.  In particular, in \cite{NP} Novruzi and Pierre proved a quite general and helpful structure result 
about second order shape derivatives.  However, this requires the a priori assumption that they exist, and demands  that the boundary of the domain $\Omega$ where the derivatives are computed is sufficiently smooth (for the 
structure of shape derivatives around irregular domains, see \cite{LaPi}). 
Moreover,  the linear and bilinear forms which appear in this structure result have to be identified  each time for the particular functional under study, and this turns out to be a delicate task (see the discussion in \cite[Section 5.9.4]{HP}).

In this work we propose a new approach to second order shape derivatives, which applies under rather mild regularity assumptions to shape functionals of the form
\begin{equation}\label{J}
J(\Omega):=-\inf\left\{ \int_{\Omega} [f(\nabla u) + g(u) ]\,dx \ :\ u\in H(\Omega)\right\}\,.
\end{equation}
Here $f$ and $g$ are assumed to be convex and smooth functions satisfying suitable growth conditions, while
$H(\Omega)$ indicates the space of functions in $H^1(\Omega)$ 
which satisfy the Dirichlet condition $u = 0$ on a fixed, nonempty, measurable portion $\Gamma _D \subset \partial \Omega$;
the Neumann part of the boundary, {\it i.e.} the complement $\partial \Omega \setminus \Gamma _D$, will be denoted by $\Gamma _N$.

As  a natural continuation of our previous paper \cite{BFL}, where we introduced a
new approach to first order shape derivatives for functionals of the type \eqref{J}, 
here we tackle second order shape derivatives by the same method.  
The main results obtained in \cite{BFL} are briefly recalled in Section \ref{secprel}, see eq.\ \eqref{J'}-\eqref{defA}. 
The basic ingredients we employ are the  analysis of the $\Gamma$-convergence  of the differential quotients appearing in the definition of $J'' (\Omega, V)$, and the duality principle 
 \begin{equation}\label{J=J*}
J(\Omega)=J^*(\Omega):=\inf\left\{ \int_\Omega [f^*(\sigma) + g^*(\div \sigma)]\,dx\ :\ \sigma \in X(\Omega;\re^n)\right\}\, ,
\end{equation}
where $X(\Omega;\re^n)$ is a suitable space of vector fields. 

Our main results, which are intimately related to each other, are:

\medskip
-- an existence result for $J'' (\Omega, V)$ which is a quadratic form in $V$ and can be  represented as follows:
\begin{equation}\label{formulaintro}
J'' (\Omega, V)=\displaystyle \int_{\partial \Omega} \big (\mathcal C(\u , V)\cdot n  \big )\, d \mathcal H ^ {n-1}
+ q(\u,V)\,; 
\end{equation}
here $\mathcal C(\u, V)\cdot n$ is the normal trace of a suitable vector field (see \eqref{defC}) depending on the solution $\u$ to $J (\Omega)$ and quadratically 
on the deformation $V$, whereas $q (\u, V)$ is a nonlocal term which involves a further vector field $ B (\u, V)$ and a quadratic form 
$ \mathcal Q (\u, \cdot) $ depending on the second derivatives of $f$ and  $g$ 
(see Theorem \ref{bdry}); we stress that, as detailed in Remark \ref{repre} below, formula \eqref{formulaintro} fits   
the general representation result for second order shape derivatives given in \cite[Corollary 2.4]{NP};

\medskip
-- a regularity result of type $W ^ {2, 2} _{\rm loc}$ for the solution to $J (\Omega)$
(see Proposition \ref{theo-reg});

\medskip
-- a new necessary optimality condition, 
which involves the distributional divergence of the above mentioned vector field  $B (\u, V)$
(see Proposition \ref{theo-div}).

\medskip
Both the regularity of the solution and the optimality condition on one hand seem to have an independent interest, and  on the other hand are strictly related to the second order shape derivative. Namely, if the $W ^ {2,2} _{\rm loc}$ regularity of $\u$ stated in Proposition \ref{theo-reg} extends up to the boundary, the field $\mathcal C (\u, V)$ turns out to admit 
a normal trace on the boundary as soon as the latter is piecewise $C ^1$; moreover, in order to arrive at the expression \eqref{formulaintro} for $J'' (\Omega, V)$, 
the optimality condition given in Proposition \ref{theo-div} is exploited as a crucial tool.

\medskip
The paper is organized as follows. After providing some notation and preliminary background in Section \ref{secprel}, we state  
in Section \ref{secmain}  our main results (Proposition \ref{theo-reg}, Proposition \ref{theo-div} and Theorem \ref{bdry}), followed by some comments and examples. In fact, we show that we can recover straightforward the second order shape derivative of the torsional rigidity, as given {\it e.g.}\ in \cite[Section 5.9.6]{HP}, and we are able to extend the formula to the case of mixed Dirichlet-Neumann conditions. 
Moreover, we consider the case of $p$-torsional rigidity for $p>2$. In this case, to the best of our knowledge, only the first order shape derivative is available in the literature (see \cite{ChMa, FaEm}). We are able to compute the second order shape derivative under an additional assumption, which is quite technical and is related to a deep regularity issue for the $p$-Laplace operator; checking its validity seems to be a delicate problem of independent interest which goes beyond the scope of this paper.

The next three sections are devoted to the proofs: the existence of the second order shape derivative (Theorem \ref{bdry} (i)) is achieved in Section \ref{secexi};
the regularity result of Proposition \ref{theo-reg} and the optimality conditions of  Proposition \ref{theo-div} are proved in Section \ref{secreg}; 
eventually, in Section \ref{secrep} we prove the representation formula for the second order shape derivative (Theorem \ref{bdry} (ii)). 

In Section \ref{secvar} we present some variants of our results, including in particular the case of the $p$-torsional rigidity for $p \in [2, + \infty)$, and then we conclude the paper by addressing some possible perspectives. 
Section \ref{secapp} is an appendix where some auxiliary lemmas and technical facts are collected. 

\bigskip
{\bf Acknowledgments.} We thank B.~Sciunzi for pointing out some useful bibliographical references. 
This work  was supported by the University of Toulon and GNAMPA (INDAM). We  gratefully acknowledge both these institutions.

\section{Preliminaries}\label{secprel}

{\it  Standing assumptions}.
Unless otherwise stated, we work under the following hypotheses, which will be referred to as standing assumptions:
\begin{itemize}
\item[--] $\Omega$ is an open bounded connected set, with a piecewise $C ^1$ boundary and unit outer normal $n$;
\item[--] $f:\re^n \to \re$ and $g: \re\to \re$ are of class $C^2$ and are strongly convex, namely there exist positive constants $m$ and $k$ such that 
\smallskip
\begin{equation}\label{km}
 (\nabla^2 f  - m I) \hbox{  is positively semidefinite, } \qquad \qquad  g'' - k \geq 0;\qquad  
\end{equation}
\item[--] $f$ and $g$ satisfy the growth conditions
\begin{equation}\label{growth}
\begin{cases}
&C_1 ( |z|^2 - 1)  \leq  f(z) \leq C_2 (|z|^2 + 1) \qquad \forall z \in \re ^n
\\ \noalign{\medskip}
& C_3 ( |v|^2 - 1)  \leq g(v) \leq C_4 (|v|^{q}+1) \qquad \  \forall v \in \re\, ,
\end{cases}
\end{equation}
where $C_i$ are positive constants, while the exponent $q$ is assumed to satisfy 
 $
 q= 2 ^* : = \frac{2n}{n-2}$  if  $n>2$ and $ q \in (1, +\infty)$   if   $n \leq 2$;

\item[--] $g(0) = 0$ (which is not restrictive up to a translation);
\item[--] the solution $\u$ to $J(\Omega)$ is Lipschitz (we point out that, by strict convexity, $\u$ is unique, and we refer to \cite{Ce, FiTr, MaTr, Mi}  for some Lipschitz regularity results);
\item[--] the deformation field $V$ belongs to $C^1(\re^n;\re^n)$.
\end{itemize}

\medskip
{\it Some standard notation.}
We adopt the convention of repeated indices. Given two vectors $a,b$ in $\re^n$
we use the notation  $\langle a , b\rangle$ to denote their Euclidean scalar product; moreover we denote by $a\otimes b$ the tensor product of $a$ and $b$, namely the matrix $(a\otimes b)_{ij}:=a_i b_j$. Given two matrices $A$ and $B$ in $\re^{n \times n}$, we write $A:B$ to denote their Euclidean scalar product, namely $A:B= A_{ij} B_{ij}$. We denote by $A^{-1}$ and $A^{T}$ the inverse and the transpose matrices of $A$, by $A^{-T}$ the transposition of the inverse of $A$ and by $I$ the identity matrix. Moreover, we denote by $a_k(A)$ the $k$-th invariant of the matrix $A$, in particular 
$$a_1(A)={\rm tr (A)} \, , \qquad  a_2 (A) = \frac{1}{2} \big ( {\rm tr} (A) I - A ^ T ) : A \, , \qquad a_n(A)=\det A\,.$$
For a tensor field $A\in C^1(\re^n;\re^{n\times n})$, by $\div A$ we mean its divergence with respect to lines, namely
$(\div A)_i:= \partial_j A_{ij}$.


Let $\Omega$ be as in the standing assumptions. 
Given a vector field $V$ on $\partial \Omega$, we decompose it into a {\it tangential and normal component} as
$V = V_\Gamma + V_n n$, with $V_n:= \langle V, n \rangle$. 

We recall that, for $\varphi \in C^1(\partial \Omega)$ and $V\in C^1(\partial \Omega;\re^n)$, the {\it tangential gradient} of $\varphi$ and the {\it tangential divergence} of $V$ are given respectively by $\nabla_\Gamma \varphi := \nabla \tilde{\varphi} - \langle \nabla \tilde{\varphi} , n \rangle n$ and $\div_\Gamma V:=\div \widetilde{V} - \langle D\widetilde{V} n ,n\rangle $ on $\partial \Omega$, where $\tilde{\varphi}$ and $\widetilde{V}$ are arbitrary $C^1$ extensions to $\re^n$ of $\varphi$ and $V$.
In particular, if $\Omega$ is of class $C^2$, the tangential divergence of $n$ gives the scalar mean curvature, that is $H_{\partial \Omega}:= \div_\Gamma n$. 
The first and second order normal derivatives of $\varphi$  will be denoted by $\partial_n \varphi := \langle \nabla \varphi , n \rangle$ and $\partial^ 2_{nn} \varphi := \langle (\nabla^2 \varphi) n , n \rangle$.

Given a vector field $\Psi \in L ^ 2 (\Omega; \re ^n)$ with distributional divergence  in $L ^ 2 (\Omega)$, 
we denote by $\Psi \cdot n$ its {\it normal trace} on $\partial \Omega$, meant as the unique element  in $H ^ {-1/2}(\partial \Omega)$ such that
\begin{equation}\label{teo-div}
\int_{\partial \Omega} (\Psi\cdot n )\,\varphi\,d{\mathcal  H^{n-1}}=\int_\Omega \big (\langle \Psi ,  \nabla \varphi \rangle  +\varphi\,\div\Psi \big ) \, dx \qquad \forall \varphi \in H ^1 (\Omega)\, ,
\end{equation}
where the boundary integral at the l.h.s.\ is intended as the duality bracket between $H ^ {-1/2} (\partial \Omega)$ and $H ^ {1/2} (\partial \Omega)$. 
For a detailed account on the theory of weak traces, we refer to \cite{A, CF} (see also \cite[Section 2.1]{BFL}).

\medskip {\it Existence of solutions, dual formulation and optimality conditions.}
Under the standing assumptions, the infimum problem $J(\Omega)$ admits a unique solution  ({\it cf.}\ \cite[Lemma 2.1]{BFL}), denoted by $\u$, in the space
$$
H(\Omega):=\{u\in H^1(\Omega)\ :\ u=0 \ \text{on }\Gamma_D\}\,.
$$
Moreover, by standard duality arguments ({\it cf.}\ \cite[Lemma 2.2]{BFL}), it holds $J(\Omega)=J^*(\Omega)$,
where $J^*(\Omega)$ is the infimum problem introduced in (\ref{J=J*}), set over the space
\begin{equation}\label{defX}
X(\Omega;\re^n):=\{\sigma \in L^2(\Omega;\re^n)\ :\ \div \sigma \in L^{q'}(\Omega)\,,\ \sigma\cdot n = 0\ \text{on }\Gamma_N\}\,,
\end{equation}
being $q' := q/(q-1)$. Also $J^*(\Omega)$ admits a unique solution, denoted by $\osigma$, which is related to $\u$ by the following differential equalities:
\begin{equation}\label{subdif}
\left\{ \begin{array}{lll}
{\osigma}= \nabla f(\nabla {\u})\ \ \hbox{a.e.\ in\ }\Omega
\\
\div {\osigma} = g'({\u})\ \ \hbox{a.e.\ in\ }\Omega\  \end{array} \right.
\quad \hbox{or equivalently}\quad
\left\{ \begin{array}{lll}
\nabla \u= \nabla f^*({\osigma})\ \ \hbox{a.e.\ in\ }\Omega
\\
\u = (g^*)'({\div {\osigma}})\ \ \hbox{a.e.\ in\ }\Omega\,.   \end{array} \right.
\end{equation}
Note that, by (\ref{subdif}) and in view of the regularity assumed on $f$, $g$ and $\u$, we have $\nabla \u, \osigma \in L^\infty (\Omega;\re^n)$ and $\u,\div \osigma\in L^\infty(\Omega)$. 

\medskip
{\it First order shape derivative.} 
Let $V\in C^1(\re^n;\re^n)$ be a deformation field and set for every $\e$
\begin{equation}\label{Psie}
\Omega_\e:=\Psi_\e(\Omega)\,,\quad \text{being } \Psi_\e(x):=x + \e V(x)\,.
\end{equation}
Under the standing assumptions, the
first order shape derivative of $J$ at $\Omega$ in direction $V$, defined by
\begin{equation}\label{defJ'}
J'(\Omega,V):=\lim_{\e\to 0 } \frac{J(\Omega_\e)-J(\Omega)}{\e }\,,
\end{equation}
exists and is given by 
\begin{equation}\label{J'}
J'(\Omega, V)= \displaystyle \int_\Omega A (\u) : DV\,dx = \int_{\partial \Omega} \langle A (\u)\, n, V \rangle \, d \mathcal H ^ {n-1}\,,
\end{equation}
where $A (\u)$ denotes the following tensor which turns out to be divergence free on $\Omega$ with a normal trace  $A (\u)\, n$  in $L^\infty(\partial \Omega)$:
\begin{equation}\label{defA}
A (\u) = \nabla  \u \otimes \osigma   - (f(\nabla \u)+g( \u))\,I
= \nabla  \u \otimes \osigma   + (f^*(\osigma)+g^*(\div \osigma) - 
\langle \nabla \u , \osigma \rangle - \u \div \osigma)\,I
\end{equation}

(with  $\osigma = \nabla f (\nabla \u)$). 
For the the proof of \eqref{J'}, see Theorem 3.3 and Theorem 3.7 in \cite{BFL}; actually, even though 
such results  are not stated for the case of mixed Dirichlet-Neumann boundary conditions, it can be easily checked that their proofs continue to work unaltered in such case.  
Moreover, if $\u$ is Lipschitz and $\nabla \u  , \, \osigma$ are in $BV(\Omega;\re^n)$, the boundary integral in \eqref{J'} can be rewritten as a linear expression of $V_n$, 
 namely it holds $J'(\Omega, V)= l_1(V_n)$ with  
$$
l_1(\varphi):=\int_{\Gamma_D} f^*(\osigma)\,  \varphi\, d\mathcal H^{n-1}    -    \int_{\Gamma_N}(f(\nabla \u) + g(\u) )\, \varphi \,d\mathcal H^{n-1}\,.
$$

\medskip
{\it Reformulation from a variable domain to a fixed domain.} Following the same procedure adopted in \cite[Lemma 4.1]{BFL} in order to recast the first order shape derivative $J'(\Omega,V)$, we rewrite the variational problem $J(\Omega _\e)$ over the fixed domain $\Omega$: by using a standard change of variables, we obtain
\begin{equation}\label{Je2}
J(\Omega_\e) = \displaystyle - \inf_{u\in H(\Omega)} \int_\Omega [f_\e(\nabla u) + g_\e(u)]\,dx
\end{equation}
where $f _\e$ and $g _\e$ are given by
\begin{equation}\label{fge} 
f_\e(x,z)  : =  f(D\Psi_\e^{-T}z)\beta_\e\, , \qquad 
g_\e(x,v)  :=   g(v)\beta_\e \, 
\end{equation}
being  $\beta_\e$ the Jacobian of the map $\Psi_\e$. Note that, in view of (\ref{Psie}),  $\beta_\e$ can  be written as  a  polynomial in $\e$, with the invariants of $DV$ as coefficients:
\begin{equation}\label{betae}
\beta _\e := |\det D\Psi_\e| = 1 + a _1 (DV) \e + a _2 ( DV) \e ^ 2 + \dots + a_n (DV) \e ^n \,.
\end{equation}

In a similar way, the dual variational problems $J^*(\Omega _\e)$ can be rewritten as 
\begin{equation}\label{J2*} 
J^*(\Omega_\e) = \displaystyle  \inf_{\sigma \in X(\Omega;\svre^n)} \int_\Omega [f_\e^*(\sigma) + g_\e^*(\div \sigma)]\,dx \, ,
\end{equation}
where $f_\e^*$ and $g_\e^*$ are the Fenchel conjugates of $f_\e$ and $g_\e$ with respect to the second variable, which by direct computation satisfy 
\begin{equation}\label{fge*}
f_\e^*(x,z)  =f^*(\beta_\e^{-1} D\Psi_\e z)\beta_\e\, , \qquad
g_\e^*(x,\xi) = g^*(\beta_\e^{-1} \div \xi)\beta_\e\, .
\end{equation}
In the subsequent asymptotic analysis as $\e \to 0$, we shall exploit the following developments holding for small $\e$, 
being $\Psi _\e$ and $\beta _\e$ defined respectively by \eqref{Psie} and \eqref{betae}: 
\begin{eqnarray}
& D\Psi_\e^{-T}=I - \e DV^T + \e^2 (DV^T)^2 + \e^3 M_\e\, , \quad \text{ with } \sup_\e \|M_\e\|\leq C\, ,
\label{Me} \\ \noalign{\medskip}
& \beta_\e =  1 + \e \div V  + \e^2 a_2(DV) + \e^3 m_\e \, ,\quad \text{ with }  \sup_\e \|m_\e\|\leq C\, ,
\label{be} \\ \noalign{\medskip}
&\beta_\e^{-1}= 1 - \e \div V + \e^2 ((\div V)^2 - a_2(DV)) + \e^3 \wt{m_\e}(V)\, ,\quad \text{ with }   \sup_\e \|\wt{m_\e}\|\leq C\, .
\label{be-1}
\end{eqnarray}

\section{Main results}\label{secmain}

All the results in this Section concern the shape functional $J$ defined in \eqref{J} and are stated under the standing assumptions given in Section \ref{secprel}. 
Our goal is to provide an existence and representation result for the second order shape derivative of $J$, meant according to the following

\begin{definition}
Given $V\in C^1(\re^n;\re^n)$,  $J$ is said to be {\it second order differentiable} at $\Omega$ in the direction $V$, if the following limit exists:
\begin{equation}\label{defJ''}
J''(\Omega,V):=\displaystyle{\lim_{\e\to 0 } 2\, \frac{J(\Omega_\e)-J(\Omega)- \e J'(\Omega, V)}{\e^2} }\,,
\end{equation}
where $J'$ is the first order shape derivative of $J$ at $\Omega$ in direction $V$  (see \eqref{defJ'}  and \eqref{J'}
for its definition and representation formula). 
\end{definition}

To prepare our main result, some preliminaries are in order. The next two  propositions will be necessary in order to prove the main theorem, but also seem  to have an autonomous interest: 
as mentioned in the Introduction, the former  concerns the regularity of the unique solution
$\u$ to $J (\Omega)$ (for $W ^ {2,2}$-type regularity results for solutions to 
quasilinear elliptic equations, see for instance \cite{DaSc, Sc, Sci} and references therein);
the latter states a new necessary condition for optimality.

\begin{proposition}\label{theo-reg} {\rm (regularity of the solution)}

\smallskip
Under the standing assumptions, $\u$ belongs to
$W^{2,2}_{\rm loc}(\Omega)$.
\end{proposition}

\medskip

\begin{proposition}\label{theo-div}  {\rm (necessary optimality condition)}

\smallskip Under the standing assumptions,
the vector field
defined by \begin{equation} \label{defB}
B (\u, V)  := \displaystyle  \nabla ^2 f(\nabla\u) (\nabla ^ 2 \u) V -  (DV-\div V\,I)\osigma
\end{equation}
satisfies \begin{equation}\label{div_eta}
\div \left[B (\u, V)\right] = g''(\u) \langle V , \nabla \u \rangle  + g'(\u) \div V \qquad in\ \mathcal D'(\Omega)\,.
\end{equation}

\end{proposition}

\bigskip
The vector field $B (\u, V)$ defined in \eqref{defB} will appear in the representation formula for the second order shape derivative. Besides $B (\u, V)$, such formula will involve another vector field depending on $\u$ and  $V$, and a quadratic form on the space $H (\Omega)$, which are defined respectively as follows. 

We let $\mathcal C (\u, V)$ be the following vector field, which  depends quadratically on  $V$:
\begin{equation}\label{defC}
\begin{array}{ll}\mathcal C (\u, V):=&\displaystyle
-\langle V , \nabla \u \rangle ( \nabla ^ 2 f (\nabla \u)) ( \nabla ^ 2 \u) \, V - \langle V , \nabla \u \rangle  (\div  V I - DV)  \nabla f(\nabla \u) \, + 
\\ \noalign{\smallskip}
& +\,  \langle DV \, \nabla f(\nabla \u) , \nabla \u \rangle V - \langle DV\, V , \nabla \u \rangle \nabla f(\nabla \u) \, +
\\ \noalign{\smallskip}
& -\, \big (  f (\nabla \u) + g (\u) \big )  (\div V I -  DV   ) V\,.
\end{array}\end{equation}
We emphasize that the distributional Hessian of $\u$ appearing in \eqref{defB} and \eqref{defC} belongs to $L ^ 2_{\rm loc}(\Omega)$ thanks to
Proposition \ref{theo-reg}. 

We let  $\mathcal Q(\u, \cdot)$ be the quadratic form given by
\begin{equation}\label{defQ}
 \mathcal Q (\u, w) :=\displaystyle  \int _{\Omega}2\,  [Q _f (x, \nabla w ) + Q _g (x, w) ] \, dx \qquad \forall w \in H (\Omega) \, ,
\end{equation}
where
$Q_{f}(x, \cdot)$  and $Q_g(x, \cdot)$ denote
the quadratic integrands
$$
Q_{f}(x,z):=\frac{1}{2}\,\langle \nabla^2 f(\nabla \u(x)) z,z \rangle\quad  \forall z \in \re ^n, \qquad  \ Q_g(x,v) := \frac{1}{2}\, g''(\u(x)) v^2\quad \forall v \in \re\,.
$$
In the sequel for brevity we shall omit to denote the dependence of $Q_f$ and of $Q_g$ on $x$ and we shall simply write
$Q_f (\cdot)$ and $Q_g (\cdot)$ in place of $Q_f (x, \cdot)$ and $Q_g(x, \cdot)$.  

\medskip

We are now in a position to state our main result.

\medskip

\begin{theorem}\label{bdry} {\rm (existence and representation of the second order shape derivative)}

\smallskip
Under the standing assumptions there hold:
\begin{itemize}
\item[(i)]  the second order shape derivative $J'' (\Omega, V)$ defined in \eqref{defJ''} exists and it is a quadratic form in $V$; 

\item[(ii)] if $\u$ is $W ^ {2,2}(\Omega)$, the vector fields  $\mathcal C (\u, V)$ and $B(\u,V)$, defined respectively in \eqref{defC} and \eqref{defB}, admit a normal trace in $H^{-1/2}(\partial\Omega)$, and
$J'' (\Omega, V)$ is given by 
\begin{equation}\label{J''1}
J'' (\Omega, V)=\displaystyle \int_{\partial \Omega} \big (\mathcal C(\u , V)\cdot n  \big )\, d \mathcal H ^ {n-1}
+ q(\u,V)\,,
\end{equation}
where $q(\u,V)$ is the following non local term:
\begin{equation}\label{defq}
q(\u,V):= - \inf _{w \in H (\Omega)}\Big \{ \mathcal Q  \big (   \u,  w - \langle V , \nabla \u \rangle \big )  + 2\int_{\partial \Omega}\big(w-\langle V , \nabla \u \rangle\big)\, B (\u, V) \cdot n \, d \mathcal H ^ {n-1} 
\Big \}\,.
\end{equation}
\end{itemize}
\end{theorem}

\begin{remark}\label{alternativeC}{\rm
Let us show how the representation formula for $J''(\Omega, V)$ given by \eqref{J''1} can be rewritten if one prefers to separate 
the contributions coming from the Dirichlet and Neumann portions of the boundary. 
First observe that, since $\u$ is assumed to be $W ^ {2,2}$ up to the boundary, the product $\langle V , \nabla \u \rangle$ belongs to $H^1(\Omega)$;  thus, by considering the traslation $v:= w - \langle V , \nabla \u \rangle $ in (\ref{defq}), the nonlocal term $q(\u,V)$ can be recast as
\begin{equation}\label{defq2}
q(\u,V)=-\inf _{v \in H^1 (\Omega)\atop v = - V_n \partial_n \u\ \text{on }\Gamma_D} \Big \{  \mathcal Q ( \u,  v) + 2\int_{\Gamma_N }v\, B (\u, V) \cdot n \, d \mathcal H ^ {n-1} 
\Big \}   + 2\int_{\Gamma_D}\langle V , \nabla \u \rangle\, B (\u, V) \cdot n \, d \mathcal H ^ {n-1}\,.
\end{equation}
By combining (\ref{J''1}) and (\ref{defq2}), we obtain
\begin{equation}\label{J''2}
\begin{array}{ll}
J''(\Omega,V) & = \displaystyle \int_{\Gamma_D} \big (\mathcal C(\u , V) + 2 \langle V , \nabla \u \rangle\, B (\u, V)\big)\cdot n \, d \mathcal H ^ {n-1} + 
\int_{\Gamma_N} \big (\mathcal C(\u , V)\cdot n \big) \, d \mathcal H ^ {n-1} + 
\\ \noalign{\smallskip}
& \displaystyle -\inf _{v \in H^1 (\Omega)\atop v = - V_n \partial_n \u\ \text{on }\Gamma_D} \Big \{  \mathcal Q ( \u ,  v) + 2\int_{\Gamma_N }v\, B (\u, V) \cdot n \, d \mathcal H ^ {n-1} 
\Big \} \,; 
\end{array}
\end{equation}
in Remark \ref{repre} below we will see that, under suitable regularity assumptions,  the infimum problem appearing in \eqref{J''2} is a quadratic form in $V_n$.

Next, in order to further simplify (\ref{J''2}), we observe that the boundary integrals appearing therein depend only on the normal traces of $\mathcal C(\u,V)$ and $B(\u,V)$ on $\partial \Omega$. 
Therefore, such fields can be replaced by any simpler vector field with the same normal trace, which can be chosen in different ways on the Dirichlet and Neumann portions of the boundary. On $\Gamma _D$, since $\nabla \u$ is parallel to $n$ on $\partial \Omega$, $\mathcal C(\u , V) + 2 \langle V , \nabla \u \rangle\, B (\u, V)$ has the same normal trace as
\begin{equation} \label{CD} 
\mathcal C_D(\u,V):=\langle V , \nabla \u \rangle ( \nabla ^ 2 f (\nabla \u)) ( \nabla ^ 2 \u) \, V + \big( \langle \nabla f(\nabla\u),  \nabla \u  \rangle -  f (\nabla \u)  \big) (\div  V I - DV) V
\, ; 
\end{equation}
while on $\Gamma_N$, since $\osigma\cdot n =0$ on such portion of the boundary, $\mathcal C(\u,V)$ has the same normal trace as 
$$\begin{array}{ll}\mathcal C_N (\u, V):=&\displaystyle
 -\langle V , \nabla \u \rangle ( \nabla ^ 2 f (\nabla \u)) ( \nabla ^ 2 \u) \, V + \langle V , \nabla \u \rangle  DV  \nabla f(\nabla \u) \, + 
\\ \noalign{\smallskip}
& + \langle DV \, \nabla f(\nabla\u), \nabla \u \rangle V -\, \big (  f (\nabla \u) + g (\u) \big )  (\div V I - DV ) V\,.
\end{array}
$$
We end up with the following reformulation of \eqref{J''1}:
\begin{equation}\label{J''3}
\begin{array}{ll}
J''(\Omega,V)&= \displaystyle \int_{\Gamma_D} \big (\mathcal C_D(\u , V)\cdot n \big)\, d \mathcal H ^ {n-1} + 
\int_{\Gamma_N} \big (\mathcal C_N(\u , V)\cdot n \big) \, d \mathcal H ^ {n-1} +
\\ \noalign{\smallskip}
& \displaystyle -\inf _{v \in H^1 (\Omega)\atop v = - V_n \partial_n \u\ \text{on }\Gamma_D} \Big \{  \mathcal Q ( \u,  v) + 2\int_{\Gamma_N }v\, B (\u, V) \cdot n \, d \mathcal H ^ {n-1} 
\Big \} \,;
\end{array}
\end{equation}
note that we readily obtain the expression of the second order shape derivative in the pure Dirichlet case, {\it i.e.} when $\Gamma _D = \partial \Omega$, by dropping the integrals over $\Gamma _N$ in the above equality. 
}
\end{remark}

\begin{remark}\label{repre}{\rm 
It is interesting to  compare our representation formula \eqref{J''1} with the structure result proved in \cite{NP}, 
which allows to enlighten the role played by the tangential and the normal components of the deformation field $V$ on $\partial\Omega$. Notice that, in general, these components  are not decoupled, {\it i.e.} $J''(\Omega, V) \not= J''(\Omega, V_\Gamma) + J''(\Omega, V_n n)$. In fact,  according to the structure theorem given in \cite[Corollary 2.4, Remark 2.10]{NP} (see also \cite[Theorem 5.9.2]{HP}), under suitable regularity conditions $J''(\Omega,V)$ can be decomposed as 
$
J''(\Omega,V)= l_1(z) + l_2(V_n)
$,
where $l_2$ is a quadratic form on $C^1(\partial \Omega)$ and $l_1$ is the linear form on $C^0(\partial \Omega)$ associated with the first order shape derivative, evaluated in 
\begin{equation}\label{defz}
z:= \langle V_\Gamma , D n \, V_\Gamma \rangle - 2 \langle  V_\Gamma , \nabla _\Gamma V_ n\rangle\,;
\end{equation}
in this way the scalar $z$ encodes the coupling between the tangential and normal components of $V$.

In the case of our functional $J$,  by Theorem 3.7 in \cite{BFL}, the linear form $l_1$ is given by 
\begin{equation}\label{defl1}
l_1(\varphi)=\int_{\Gamma_D} f^*(\osigma)\,  \varphi\, d\mathcal H^{n-1}    -    \int_{\Gamma_N}(f(\nabla \u) + g(\u) )\, \varphi \,d\mathcal H^{n-1}\,,
\end{equation}
and we claim that, as soon as $\partial\Omega$ is  $C ^2$,  by Theorem \ref{bdry} the quadratic form $l _2$ is given by
\begin{equation}\label{l2D}
\begin{array}{lll} l_2(\varphi) & : =  \displaystyle \int_{\Gamma_D} \varphi^2 \left[  \partial_n \u \langle D\osigma n, n \rangle + f^*(\osigma) H_{\partial \Omega}\right] \,d\mathcal H^{n-1} +
\\ \noalign{\smallskip}
& \displaystyle +
 \int_{\Gamma_N} \Big[  \langle \osigma, \nabla_\Gamma (\varphi^2 \partial_n\u)\rangle  - \varphi^2 \big[ \langle (\nabla^2\u)\,\osigma,n\rangle 
+ (f(\nabla\u) + g(\u) ) H_{\partial \Omega}
+ \partial_n\u \langle D\osigma n, n\rangle  \big]
\Big]\,d\mathcal H^{n-1}+
\\ \noalign{\smallskip}
& \displaystyle - \inf _{v \in H ^1(\Omega)\atop v = - \varphi \partial_n \u\ \text{on }\Gamma_D} \Big \{  \mathcal Q ( \u,  v) + 2\int_{\Gamma_N }v\, 
[\varphi \langle D\osigma n,n \rangle - \langle \osigma , \nabla_\Gamma \varphi \rangle]
\, d \mathcal H ^ {n-1} 
\Big \}\,.
\end{array}
\end{equation}

This claim can be easily checked when $V_\Gamma = 0$  and $\partial\Omega=\Gamma_D$.  Indeed in this case we have $z= 0 $ and hence $l_1(z)=0$; moreover, starting from equality \eqref{J''1}, if we replace  $\mathcal C(\u,V)$ by the simpler field $\mathcal C_D(\u,V)$ given in \eqref{CD} and we use the identity $ \langle(\div  V I - DV) n , n \rangle = V_n H _{\partial \Omega}$ (holding thanks to the $C^2$ regularity assumption on $\partial \Omega$), we are led to the equality $J''(\Omega,V)= l_2(V_ n)$, with $l_2$ defined according to \eqref{l2D}.
The proof of the claim in the general case requires some tedious but straightforward computations, 
that we omit for the sake of conciseness.

}
\end{remark}

\begin{remark}
{\rm  We point out that the nonlocal part of the second variation, namely the term $q(\u,V)$ defined in (\ref{defq}), may also be reformulated in dual form. 
In analogy  with (\ref{defQ}), let us introduce the quadratic form $\mathcal Q ^* (\osigma, \cdot) $ on the space  $X(\Omega;\re^n)$ defined in (\ref{defX}): we set
$$
\mathcal Q ^* ( \osigma, \eta ):=\displaystyle \int _{\Omega}2[ Q_{f^*} (x,\eta ) + Q _{g^*} (x,\div \eta ) ]\, dx\,,$$
with
$Q_{f^*}(x,z):=\frac{1}{2}\, \langle \nabla^2 f^*(\overline \sigma(x)) z,z \rangle$ and $Q_{g^*}(x,v) := \frac{1}{2}\, (g^*)''(\div \osigma(x)) v^2$. (For brevity in the sequel we adopt the notation 
$Q_{f^*}(\cdot)$ and $Q_{g^*}(\cdot)$ in place of $Q_{f^*}(x,\cdot)$ and $Q_{g^*}(x,\cdot)$.)

A standard duality argument (see Lemma \ref{prop_duality}) yields
$$
q(\u,V) = \inf_{\eta \in X(\Omega;\re^n)} \Big \{  \mathcal Q^* ( \osigma ,  \eta-B(\u,V) ) + 
\int_{\Gamma_D} 2\,  \partial _n \u \, V_n \, (\eta \cdot n) \, d \mathcal H ^ {n-1}  \Big \} \,,
$$
or equivalently, via a translation,
$$
q(\u,V)=\inf_{\eta \in X(\Omega;\svre^n)-B(\u,V)} 
\Big \{  \mathcal Q^* ( \osigma,  \eta) + 
\int_{\Gamma_D} 2\,  \partial _n \u \, V_n \, (\eta \cdot n) \, d \mathcal H ^ {n-1}  \Big \}  + \int_{\Gamma_D} 2\,  \partial _n \u \, V_n \, (B(\u,V) \cdot n) \, d \mathcal H ^ {n-1} \,.
$$
}
\end{remark}

\begin{example}
{\rm  An example of shape functional which is covered by Theorem \ref{bdry} 
is the torsional rigidity under mixed Dirichlet-Neumann boundary conditions:
$$
J(\Omega) = -\inf\left\{ \int_\Omega \Big ( \frac{|\nabla u|^2}{2} - \lambda u \Big ) \,dx\ :\ u\in H (\Omega)\right\}\, .
$$
In this case, assuming $\partial \Omega$ of class $C^2$ and taking deformations normal to the boundary, the second order shape derivative given in (\ref{J''3}) reads
\begin{equation}\label{J''torsion}
\begin{array}{lll}  
J'' (\Omega, V )  = & \displaystyle - \frac{1}{2} \int_{\Gamma_D}  V_n ^ 2 \big ( 2 \lambda \partial_n \u + |\partial_n\u |^2 H _{\partial \Omega} \big ) \, d \mathcal H ^ {n-1} 
+
\\ \noalign{\smallskip}
 & \displaystyle 
- \frac{1}{2} \int_{\Gamma_N}  V_n ^ 2 \big [ ( - 2 \lambda \u + |\nabla \u |^2 ) H _{\partial \Omega} + 2 \langle \nabla^2 \u \nabla \u , n \rangle \big ] \, d \mathcal H ^ {n-1} +
\\ \noalign{\smallskip}
& \displaystyle
- \inf_{ v \in H ^ 1 (\Omega)\atop v= - V_n\,\partial_n \u \text { on } \Gamma_D } \left\{
 \int _{\Omega} |\nabla v| ^ 2\,dx  - 2 \int_{\Gamma_N} v \langle \nabla \u , \nabla_\Gamma V_n \rangle \, d \mathcal H ^ {n-1}   \right\}\,.
\end{array}
\end{equation}
The above equality can be established by applying Theorem \ref{bdry} (as it continues to hold for $g ( t) =- \lambda t$, {\it cf.} the variant discussed in \S \ref{subsec-linear}). Namely in this case, as already pointed out in Remark \ref{repre}, the second order shape derivative is given by $l_2(V_n)$, being $l _2$ the quadratic form defined in (\ref{l2D}). The equality between $l_2(V_n)$ and the r.h.s.\ of (\ref{J''torsion}) readily follows by using the identity $\Delta \u = \partial^2_{nn}\u + \partial_n \u H _{\partial \Omega} $ on $\Gamma_D$, the fact that $\partial_n\u=0$ on $\Gamma_N$ and by taking into account that
$Q_g \equiv 0$ and $Q_f (z) = \frac{|z| ^ 2}{2}$.

We underline that in the pure Dirichlet case, {\it i.e.}\ when $\Gamma_N=\emptyset$, our formula \eqref{J''torsion} allows to recover the known expression of the second order shape derivative of the torsional rigidity
(see for instance \cite[Section 5.9.6]{HP}).
Finally let us point out that, when the deformation field is not normal to the boundary, the expression of $J'' (\Omega, V)$ is slightly more complicated than (\ref{J''torsion}). Precisely, it is given by $l_1(z) + l_2(V_n)$, being $l_1$ and $z$ defined respectively in (\ref{defl1}) and (\ref{defz}); in other words, we have to add to the r.h.s.\ of formula (\ref{J''torsion}) the following terms:
$$
\int_{\Gamma_D} \frac{|\partial_n \u|^2 }{2} (\langle V_\Gamma , D n \, V_\Gamma \rangle - 2 \langle  V_\Gamma , \nabla _\Gamma V_ n\rangle) \, d \mathcal H ^ {n-1} - \int_{\Gamma_N} \Big ( \frac{|\nabla\u|^2 }{2}-\lambda \u\Big ) (\langle V_\Gamma , D n \, V_\Gamma \rangle - 2 \langle  V_\Gamma , \nabla _\Gamma V_ n\rangle)\, d \mathcal H ^ {n-1}\, .
$$
}

\end{example}

\begin{example}\label{examplep}
{\rm 
As a variant of the previous example, 
we can consider the $p$-torsion problem under Dirichlet boundary conditions: 
\begin{equation}\label{def-p}
J_p(\Omega) = -\inf\left\{ \int_\Omega \Big ( \frac{|\nabla u|^p}{p} - \lambda u \Big ) \,dx\ :\ u\in W^{1,p}_0(\Omega)\right\}\, ,
\end{equation}
where $p> 2$ and $\lambda$ is a positive parameter.  Note that in this case the integrands $f$ and $g$ do not satisfy the standing assumptions, in particular (\ref{km}), so that we cannot apply directly 
Theorem \ref{bdry}. In fact a major difficulty appears with the degeneracy at the origin of the Hessian of $f$. 
Nevertheless, by exploiting a suitable approximation argument, we are able to 
 show that Theorem \ref{bdry} is still valid,  provided a suitable equality between weighted Sobolev spaces holds ({\it cf.}\ eq.\ \eqref{HW}).  The outcoming formula reads as follows. Let $\partial \Omega$ be of class $C^2$ and let the deformation field $V$ be normal to the boundary.  
Observe that $\osigma$ is parallel to $\nabla \u$, and satisfies $\div \osigma = -\lambda$; moreover, the equalities
$\langle (\div \osigma I - D\osigma ) n , n \rangle = \osigma_n H_{\partial \Omega}$
and $\Delta u = \partial^2_{nn} \u + \partial_n \u H _{\partial \Omega} $ hold on $\partial \Omega$. 
Taking these facts into account, starting from \eqref{J''3} we obtain: 
\begin{equation}\label{formulap}
 J''_p (\Omega, V)=- \frac{1}{p} \int_{\partial \Omega}  V_n ^ 2 \big ( p \lambda \partial_n \u  + |\partial_n \u| ^ p H _{\partial \Omega} \big ) \, d \mathcal H ^ {n-1}  
- \inf_{ v \in H^{1} (\Omega)\atop v= - V_n \partial _n \u\text { on } \partial \Omega }
 \int _\Omega
\langle P (\u)  \nabla v , \nabla v \rangle\,dx\,,
\end{equation}
with 
$$P(\u) :=  |\nabla\u|^{p-2}\left( I + (p-2) \frac{\nabla \u}{|\nabla\u|} \otimes \frac{\nabla \u}{|\nabla\u|}\right)\,. \
$$
The complete presentation of this variant is postponed to \S \ref{sec-p}, see Theorem \ref{prop-p}. 
}
\end{example}

\section{Existence of the second order shape derivative}\label{secexi}

This section is devoted to the proof of Theorem \ref{bdry} (i).
We introduce for brevity the differential quotients 
\begin{equation}\label{r_e}
r_\e(V):=\displaystyle{ 2\, \frac{J(\Omega_\e)-J(\Omega)- \e J'(\Omega, V)}{\e^2} }\,.
\end{equation}
By exploiting the formulation (\ref{Je2})  of $J (\Omega _\e)$ and taking therein as a test function $\u + \e w$, we infer
\begin{equation}\label{initial_lb}
 r_\e(V) \geq \frac{2}{\e^2} \left(  - \int_\Omega [f_\e (\nabla \u + \e \nabla w) + g_\e (\u + \e w) ]\,dx - J(\Omega,V) -\e J'(\Omega,V) \right)\, .
\end{equation}

We are thus led to introduce the sequence of functionals defined for $w \in H (\Omega)$  and  $V\in C^1(\re^n;\re^n)$ by
\begin{align}
E_\e(w, V) & :=\frac{2}{\e^2} \left\{\int_\Omega [f_\e (\nabla \u + \e \nabla w ) + g_\e (\u + \e w)]\, dx + J(\Omega) + \e J'(\Omega, V)\right\} \,. \label{E_e}
\end{align}
Moreover,  still for $w \in H (\Omega)$ and  $V\in C^1(\re^n;\re^n)$,
 let us define the following functional $E (w, V)$, which will be seen to be the variational limit of $E _\e(w, V)$, and its ``dual counterpart'' $E ^*(\eta, V)$
(with $\eta \in X(\Omega;\re^n)$ and $V$ as above):
\begin{equation}\label{E}
\begin{array}{ll}
E(w, V) & := 
 \displaystyle 2  \int_\Omega [f(\nabla \u) + g(\u)]a_2(DV) \,dx + 2 \left\{ \int_\Omega \left [ Q_f(\nabla w - DV^T \nabla \u) + Q_g(w)\right ] \,dx \,+ \right.
\\ \noalign{\medskip}
& \displaystyle \left . + \int_\Omega \left[  - \langle (DV-\div V\,I)\nabla f(\nabla \u) , \nabla w - DV^T \nabla 	\u\rangle\,dx  
\,+  \div V g'(\u) w\right]  \, dx \right\} \, ,
\end{array}
\end{equation}
\begin{equation}\label{E*}
\begin{array}{ll}                   
E^*(\eta, V) 
 &:=  \displaystyle2 \int_\Omega\big  [ f^*(\osigma) + g^*(\div\osigma)- \langle \nabla f^*(\osigma),\osigma \rangle -
\u \div\osigma\big ] \, a_2(DV) \,dx
+ 
\\ \noalign{\medskip}
& \displaystyle + 2  \Big\{\int_\Omega \Big [Q_{f^*}((DV - \div V \, I )\overline \sigma  + \eta) + Q_{g^*} (\div \eta - \div V \div \osigma)
  + \langle \nabla f^*(\osigma), DV\,\eta \rangle \Big ]\,dx \Big\}\,.
\end{array}
\end{equation}

Theorem \ref{bdry} (i) holds in view of the following result:

\begin{proposition}\label{thm_J''}
Under the standing assumptions, the second order shape derivative $J'' (\Omega, V)$ exists, and is given by
\begin{equation}\label{J''_primal}
J''(\Omega, V) =  \lim_{\e \to 0} \big ( - \inf_{w \in H (\Omega)} E_\e(w, V) \big ) =  -\inf_{w \in H(\Omega)} E(w, V)
= \inf_{\eta \in X(\Omega;\re^n)} E^*(\eta, V)\,.
\end{equation}
Moreover the map $ V \mapsto \, J''(\Omega, V) $ is a quadratic form in $C^1(\re^n;\re^n)$. 
\end{proposition}

\bigskip
Proposition \ref{thm_J''} will be obtained as a consequence of the two lemmas stated hereafter.

\begin{lemma}\label{propbounds} {\rm (upper and lower bounds)}

\smallskip

Under the standing assumptions, the differential quotients defined in \eqref{r_e}  satisfy
\begin{equation}\label{thesis_lb}
\liminf _{\e \to 0 } r _\e (V) \geq  -\inf_{w \in H(\Omega)} E(w, V)
\end{equation}
and
\begin{equation}\label{thesis_ub}
\limsup _{\e \to 0 } r _\e (V) \leq   \inf_{\eta \in X(\Omega;\re^n)} E^*(\eta, V)\,.
\end{equation}
\end{lemma}

\begin{lemma}\label{propgamma} {\rm ($\Gamma$-convergence)}

\smallskip

Under the standing assumptions,
for every $V \in C ^ 1 (\re ^n; \re ^n)$
the sequence $\{E_\e(\cdot, V)\}$ is equicoervice and $\Gamma$-converges to $E(\cdot, V)$, with respect to the weak topology of $H(\Omega)$. \end{lemma}

Let us show how the proposition follows from these lemmas, and then turn back to their proofs.

\bigskip
{\it Proof of Proposition \ref{thm_J''}}.
We begin by proving the last equality in \eqref{J''_primal}.
Let $\Psi: L^2(\Omega;\re^n) \times L^q(\Omega) \to \re$ be the function
\begin{equation}\label{defPsi}
\Psi(z,v):=\int_\Omega \left[\,  Q_f(z-a) + \langle z-a , b\rangle  + Q_g(v) + v \alpha + \gamma \, \right]\, dx\,
\end{equation}
with
$$a :=  DV^T \nabla \u\,,  \quad b :=  - (DV-\div V\,I)\osigma\, , \quad \alpha  :=  \div V g'(\u )\, , \quad \gamma  := [ f(\nabla \u) + g(\u) ]\, a_2(DV)\,.$$
For every $w\in H(\Omega)$, it holds $E(w,V)=2\, \Psi(\nabla w , w)$, in particular
\begin{equation}\label{eq1}
 -\inf_{w \in H(\Omega)} E(w, V) = - 2 \inf_{w\in H (\Omega)} \left\{\Psi (\nabla w ,  w)  \right\}\, .
\end{equation}
By applying Lemma \ref{prop_duality} with $Y=H(\Omega)$, $Z= L^2(\Omega;\re^n) \times L^q(\Omega) $, $A: Y \to Z$ the operator $A v:=(\nabla v, v)$, $\Phi:Y\to \re$ the zero function and $\Psi$  defined as in (\ref{defPsi}), we can rewrite \eqref{eq1} as
\begin{equation}\label{eq2}
- 2\, \inf_{w\in Y} \left\{\Psi (A w) + \Phi(w) \right\} = 2\,  \inf_{(\eta,\tau) \in Z^*} \left\{ \Psi^*(\eta,\tau) + \Phi^*(-A^* (\eta,\tau))  \right\}\, .
\end{equation}
Let us compute the Fenchel conjugates appearing in the r.h.s.\ of (\ref{eq2}).
By exploiting Lemma \ref{lemma_duality}, it is easy to see that
$$
\Psi^*(\eta,\tau) =  \int_\Omega \left[  (Q_{f})^*(\eta -b) + \langle \eta , a\rangle +   (Q_g)^*(\tau - \alpha) +\gamma \right]\, dx\, .
$$
Since $\Phi\equiv 0$, its Fenchel conjugate $\Phi^*$ is 0 at 0 and $+\infty$ otherwise.
As an element of $Y^*$, $A^*(\eta,\tau)$ is characterized by its action on the elements of $Y$: since
$$
\langle A^*(\eta,\tau),v\rangle_{Y^*,Y} = \langle (\eta,\tau), A v\rangle_{Z^* , Z} = \int_\Omega \left( 
\langle \eta ,  \nabla v  \rangle + \tau v  \right) \, dx\, ,
$$
we infer that $\Phi^*( A^*(\eta,\tau))=0$ if and only if $\tau=\div \eta$ (with the additional condition $\eta\cdot n=0$ 
 on $\Gamma_N$); note that, as a consequence, such vector fields $\eta$ belong to $ X(\Omega;\re^n)$.

Therefore we may rewrite the r.h.s.\ of (\ref{eq2}) as
\begin{equation}\label{eq3}
2\,  \inf_{\eta\in X(\Omega;\svre^n)} \left\{ \int_\Omega \left[(Q_{f})^*(\eta -b) + \langle \eta , a\rangle +   (Q_g)^*(\div\eta - \alpha) +\gamma \right] \, dx \right\}\,.
\end{equation}
Since by Lemmas \ref{lemma_quadratic} and \ref{lemma_hessian} we have
$
\left(Q_f\right)^*=Q_{f^*}$ and $\left(Q_g\right)^*=Q_{g^*}$,
and by the optimality conditions (\ref{subdif}) there holds
$$
\gamma=   \int_\Omega\big  [ f^*(\osigma) + g^*(\div\osigma)- \langle \nabla f^*(\osigma),\osigma \rangle -
\u \div\osigma\big ] \, a_2(DV) \,dx \,,
$$
we conclude that the expression in \eqref{eq3} agrees with  $\inf_{\eta \in X(\Omega;\svre^n)} E^*(\eta,V)$, which achieves the proof of the last equality in \eqref{J''_primal}. 
Combined with Lemma \ref{propbounds}, such equality implies the existence of the second order shape derivative,  and the fact that it  agrees with any of the two infimum problems at the r.h.s. of (\ref{thesis_lb}) and (\ref{thesis_ub}). 
The second equality in \eqref{J''_primal} is a direct consequence of Lemma \ref{propgamma}.

Finally, the fact that the map $V \mapsto J'' (\Omega, V)$ is a quadratic form is readily checked: one may start from the equality 
$J'' (\Omega, V) = - \inf_w E (w, V)$, and look at the expression (\ref{E}) of $E (w, V)$,  by recalling in particular that
$Q_f$ is a quadratic form.  
\qed

\bigskip
{\it Proof of Lemma \ref{propbounds}}.
We first prove the lower bound in \eqref{thesis_lb}. We start from the inequality  (\ref{initial_lb}), where we choose $w$ as an element of $C ^ \infty (\overline \Omega)\cap H (\Omega)$. 
Recalling the definitions (\ref{fge}) of $f_\e$ and $g_\e$,
the definition  (\ref{J}) of $J(\Omega)$, and the expression (\ref{J'}) of $J'(\Omega;V)$, we can rewrite the inequality (\ref{initial_lb}) separating the terms of different orders in $\e$. Recalling that $M _\e$ and $m _\e$ are defined respectively according to \eqref{Me} and \eqref{be}, and setting for brevity
$$h_1:=  \nabla w - DV^T \nabla \u\,,
\quad 
h_2:=  (DV^T)^2 \nabla \u - DV^T \nabla w\,,
\quad 
h_\e :=  (DV^T)^2\nabla w + M_\e (\nabla \u + \e \nabla w)\, ,
$$
we get
\begin{equation}\label{Ii}
 r_\e(V) \geq \Big (I_0(\e) + I_1(\e) + I_2(\e)\Big )\, ,
\end{equation}
where
\begin{align*}
I_0(\e):= & - 2 \int_\Omega [f(\nabla \u +  \e h_1 + \e^2 h_2 + \e^3 h_\e ) + g(\u + \e w)] (a_2(DV) + \e m_\e)\,dx \, , 
\\
I_1(\e):= & -  \frac{2}{\e} \int_\Omega[  f(\nabla \u +  \e h_1 + \e^2 h_2 + \e^3 h_\e ) + g(\u + \e w)  - f(\nabla \u) - g(\u) ]\div V\,dx \, , 
\\
I_2(\e):= & - \frac{2}{\e^2} \int_\Omega [ f(\nabla \u +  \e h_1 + \e^2 h_2 + \e^3 h_\e ) + g(\u + \e w) - f(\nabla \u) - g(\u) + \e \langle \nabla f(\nabla \u), DV^T \nabla \u\rangle ]\, dx                          \,.
\end{align*}

Let us study separately the asymptotic behavior of $I_i (\e)$ as $\e\to 0$.
By exploiting the growth assumptions on the integrands and the uniform $L^\infty$ boundedness of $m_\e$ and $h_\e$, we get
\begin{equation}\label{I0}
\lim_{\e\to 0} I_0 (\e) = - 2 \int_\Omega [f(\nabla \u) + g(\u)]a_2(DV) \,dx\, .
\end{equation}

By applying  Lemma \ref{lemma_taylor} (i) to the integral functionals $I_f$ and $I_g$ defined according to \eqref{Iphi} (note that this can be done after possibly modifying $g$ outside an interval, since
$\u + \e w$ remains uniformly bounded in $L^\infty$), we get
\begin{equation}\label{I1}
\lim_{\e\to 0} I_1(\e)
=
- 2 \int_\Omega [\langle \nabla f(\nabla \u), \nabla w - DV^T \nabla \u \rangle + g'(\u) w ]\div V \,dx\,.\end{equation}
Let us now consider $I_2(\e)$. We recall that, in view of the optimality conditions in \eqref{subdif}, there holds $\osigma=\nabla f(\nabla \u)$ and $\div \osigma=g'(\u)$, so that
\begin{equation}\label{addendum}
\int_\Omega [\langle \nabla f(\nabla \u) , \nabla w\rangle + g'(\u) w]\,dx  = \int_\Omega \big ( \langle \osigma ,  \nabla w \rangle + w \div \osigma \big ) \,dx = \int_{\partial \Omega} w\osigma \cdot n \, d \mathcal H ^ {n-1} = 0 \,,
\end{equation}
where the last equality holds since $w= 0$ on $\Gamma _D$ and 
$\osigma \cdot n= 0$ on $\Gamma_N$.
Thus, $I_2(\e)$ remains unchanged if we add the zero term (\ref{addendum}) (multiplied by $ 2/\e$):
\begin{align*}
I_2(\e) = & - \frac{2}{\e^2} \int_\Omega [ f(\nabla \u +  \e h_1 + \e^2 h_2 + \e^3 h_\e ) - f(\nabla \u) - \e  \langle \nabla f(\nabla u) , h_1 \rangle ]\,dx \,+ \notag \\
& - \frac{2}{\e^2}  \int_\Omega [g(\u + \e w)  - g(\u)- \e g'(\u) w]\,dx\, .
\end{align*}
Hence, by applying Lemma \ref{lemma_taylor}, we obtain
\begin{align}
\lim_{\e\to 0} I_2(\e) & = - 2 \int_\Omega [Q_f(h_1) + \langle \nabla f(\nabla \u), h_2\rangle + Q_g(w)  ]\,dx\, . \label{I2}
\end{align}
By combining (\ref{Ii}), (\ref{I0}), (\ref{I1}) and (\ref{I2}), by the arbitrariness of $w\in C ^ \infty (\overline \Omega) \cap H (\Omega)$ we infer
\begin{equation}\label{thesis_lbtilde}
\liminf _{\e \to 0 } r _\e (V) \geq  -\inf_{w \in C ^ \infty (\overline \Omega) \cap H (\Omega)} E(w, V)\,.
\end{equation}

In order to conclude the proof of the inequality \eqref{thesis_lb}, it is enough to observe that the infimum at the r.h.s.\ of (\ref{thesis_lbtilde}) coincides with the infimum at the r.h.s.\ of (\ref{thesis_lb}).  
Indeed, the integral functional $E(\cdot,V)$ above is continuous on $H(\Omega)$
(by the growth conditions \eqref{growth} and the Sobolev embedding of $W^{1,p}(\Omega)$ into  $L^p(\Omega)$), whereas
$C ^ \infty (\overline \Omega) \cap H (\Omega)$ is dense in $H(\Omega)$.  

\medskip
Let us now prove the upper bound in \eqref{thesis_ub}.
Let $\eta$ be an arbitrary element of the space $C^ \infty (\overline \Omega;\re^n) \cap X (\Omega; \re ^n)$.
In view of (\ref{J=J*}), the differential quotients $r_\e(V)$ introduced in \eqref{r_e} can be rewritten as
$$r_\e(V)=\frac{2}{\e^2}\,\left[J^*(\Omega_\e) - J^*(\Omega) - \e J'(\Omega, V) \right]\, .
$$
Then, by exploiting 
the expression of $J^*(\Omega_\e)$ in (\ref{J2*}), we infer that
\begin{equation}\label{initial_ub1}
r_\e(V) \leq \frac{2}{\e^2} \left(  \int_\Omega [f^*_\e (\osigma + \e \eta) + g^*_\e (\div \osigma + \e \div \eta) ]\,dx - J^*(\Omega) -\e J'(\Omega,V) \right)\, .
\end{equation}

Recalling the definitions (\ref{fge*}) of $f^*_\e$ and $g^*_\e$, the definition of $J ^ *(\Omega)$ in (\ref{J=J*}), and the second expression for $J'(\Omega;V)$ in (\ref{J'}), we can rewrite the inequality (\ref{initial_ub1}) separating the terms of different orders in $\e$. Recalling the definition \eqref{be-1} of $\widetilde{m_\e}$, and setting for brevity
\begin{align*}
z_1:= & ( DV - \div V I )\osigma + \eta \,,
\\
z_2:= & (DV-\div V\,I)\eta + ((\div V)^2 - \div V \,DV - a_2(DV))\osigma\,,
\\
z_\e := & \widetilde{m_\e} (\osigma + \e \eta) + ( (\div V)^2 - a_2(DV) )(DV\osigma + \eta) - \div V \, DV\eta \,,
\\
\tau_1:= &\div \eta - \div V \div \osigma\,,
\\
\tau_2:= & ((\div V)^2 - a_2(DV))\div \osigma - \div V \div \eta\,,
\\
\tau_\e := & ((\div V)^2 - a_2(DV) )\div \eta + \widetilde{m_\e} (\div\osigma + \e \div\eta)\,,
\end{align*}
we obtain
\begin{equation}\label{Ii*}
 r_\e(V) \leq \Big (I^*_0(\e)  + I^*_1(\e)  +  I^*_2(\e)  \Big )\, ,
\end{equation}
where
\begin{align*}
I_0^*(\e) := & 2 \int_\Omega [f^*(\osigma +  \e z_1 + \e^2 z_2 + \e^3 z_\e ) + g^*(\osigma+ \e \tau_1 + \e^2 \tau_2 + \e^3 \tau_\e)](a_2(DV) + \e m_\e) \,dx \, , 
\\
I_1^*(\e) := & \frac{2}{\e} \int_\Omega[  f^*(\osigma +  \e z_1 + \e^2 z_2 + \e^3 z_\e) + g^*(\osigma+ \e \tau_1 + \e^2 \tau_2 + \e^3 \tau_\e)  - f^*(\osigma) - g^*(\div \osigma) ]\div V\,dx \, , 
\\
I_2^*(\e) := & \frac{2}{\e^2} \int_\Omega [
f^*(\osigma +  \e z_1 + \e^2 z_2 + \e^3 z_\e ) + g^*(\osigma+ \e \tau_1 + \e^2 \tau_2 + \e^3 \tau_\e) - f^*(\osigma) - g^*(\div \osigma) \, +
\\ & \ \ \  \ \ - \e \langle \nabla f^*(\osigma), (DV- \div V) \osigma \rangle + \e \div V \div \osigma (g^*)'(\div \osigma)]\,dx                          \,.
\end{align*}
Let us study separately the asymptotic behavior of $I_i^*(\e)$ as $\e\to 0$.
By exploiting the growth properties of the integrands (note that the growth assumptions on $f$ and $g$ made in \eqref {growth} imply similar conditions on $f^*$ and $g^*$) and the uniform $L^\infty$ boundedness of $\widetilde{m_\e}$, we get
\begin{equation}\label{I0*}
 \lim_{\e\to 0} I_0^*(\e) =  2 \int_\Omega [f^*(\osigma )  + g^*(\div \osigma)]a_2(DV) \,dx   \,   .
\end{equation}

By applying  Lemma \ref{lemma_taylor} (i) to the integral functionals $I_{f^*}$ and $I_{g^*}$ defined according to \eqref{Iphi} (note that this can be done after possibly modifying $g^*$ outside an interval, since
$\div \osigma + \e \div \eta$ remains uniformly bounded in $L^\infty$), we get
\begin{equation}\label{I1*}
\lim_{\e\to 0} I_1^*(\e) =  2 \int_\Omega [\langle \nabla f^*(\osigma), z_1 \rangle + (g^*)'(\div \osigma) \tau_1 ]\div V   \,   dx   \,   .
\end{equation}
Let us now consider $I_2^*(\e)$. We recall that, in view of the optimality conditions \eqref{subdif},
we have
\begin{equation}\label{addendum*}
\int_\Omega [\langle \nabla f^*(\osigma ) , \eta \rangle + (g^*)'(\div \osigma) \div \eta]\,dx  = 0 \,.
\end{equation}
Thus, $I_2^*(\e)$ remains unchanged by adding the zero term (\ref{addendum*}) (multiplied by $ - 2/\e$):
\begin{align*}
I_2^*(\e) = & \frac{2}{\e^2} \int_\Omega [ f^*(\osigma +  \e z_1 + \e^2 z_2 + \e^3 z_\e ) - f^*(\osigma)
 - \e \langle \nabla f^*(\osigma), z_1 \rangle +
\\
& \ \ \  \ \ + g^*(\osigma+ \e \tau_1 + \e^2 \tau_2 + \e^3 \tau_\e) - g^*(\div \osigma)
 - \e (g^*)'(\div \osigma)\tau_1]\,dx \,.
\end{align*}
Hence, by applying Lemma \ref{lemma_taylor}, we infer
\begin{equation}\label{I2*}
\lim_{\e\to 0} I_2^*(\e)  =  2 \int_\Omega [Q_{f^*}(z_1) + \langle \nabla f^*(\osigma), z_2\rangle + Q_{g^*}(\tau_1) +  (g^*)'(\div \osigma)\tau_2 ]\,dx \,.
\end{equation}
By combining (\ref{Ii*}) with (\ref{I0*}), (\ref{I1*}) and (\ref{I2*}) and recalling the definition (\ref{E*}) of $E^*$, in view of the arbitrariness of $\eta \in C^ \infty (\overline \Omega;\re^n) \cap X (\Omega; \re ^n)$, we infer
\begin{equation} \label{final_ub1}
\limsup_{\e\to 0}  r_\e(V)   \leq \inf_{\eta \in C^ \infty (\overline \Omega;\re^n) \cap X (\Omega; \re ^n)} E^*(\eta, V)   \,   .
\end{equation}
Finally we observe that, for every $\eta \in {X} (\Omega;\re^n)$, there exists a sequence $\{ \eta _h \} \subset 
C^ \infty (\overline \Omega;\re^n) \cap X (\Omega; \re ^n)$,  such that $\eta _h \to \eta$ in $L ^ 2(\Omega;\re^n)$ and $\div \eta _h \to \div \eta$ in $L ^{q'} (\Omega)$. Since $E ^ * (\cdot, V)$ is continuous with respect to such convergence (thanks to the fact that  the growth assumptions on $f$ and $g$ made in \eqref {growth} imply similar conditions on $f^*$ and $g^*$), 
the infimum at the r.h.s.\ of (\ref{final_ub1}) turns out to agree with the infimum at the r.h.s.\  of (\ref{thesis_ub}), and    
the proof is achieved.

\qed

\bigskip
{\it Proof of Lemma \ref{propgamma}.}
For brevity, throughout the proof we fix the vector field $V$ and we omit it when writing the functionals $E_\e$ and $E$. We begin by showing the equicoercivity of $E_\e$. In view of the convexity assumptions on $f$ and $g$ we infer that the following bounds hold true, where $m$ and $k$ are the positive constants appearing in \eqref{km}: 
$$
f(z_0+z)  \geq f(z_0) + \langle \nabla f(z_0) ,  z \rangle  + m |z|^2\,,
\qquad
g(v_0 + v)  \geq g(v_0) + g'(v_0) v + k |v|^2 \,.
$$
By exploiting such bounds for $z_0=\nabla \u$, $z= - \nabla \u + D\Psi_\e^{-T}(\nabla \u + \nabla w)  $,  $v_0=\u$ and $v=\e w$, we infer that for every $\e$
$$
E_\e(w)\geq c_1 \|\nabla w\|^2_{L^2}  + c_2 \| w\|^2_{L^2} 
- c_3 \|\nabla w\|_{L^2}  - c_4  \| w\|_{L^2} - c_5\,,
$$
for some positive constants $c_i$. Therefore the sequence $E_\e$ is equicoercive in $H(\Omega)$.

\smallskip

Let us now show the $\Gamma$-convergence statement. By definition of $\Gamma$-convergence, we have to prove the so-called $\Gamma$-liminf and $\Gamma$-limsup inequalities:
\begin{eqnarray}
\inf \left\{\liminf E_\e(w_\e)\ : \ v_\e \stackrel{H(\Omega)}{\rightharpoonup} v \right \} \geq E(w)\label{G_liminf}
\\
\inf \left\{\limsup E_\e(w_\e)\ : \ v_\e \stackrel{H(\Omega)}{\rightharpoonup} v \right \}  \leq E(w)\,.  \label{G_limsup}
\end{eqnarray}

Let us prove (\ref{G_liminf}). For every $w\in H(\Omega)$ and for every sequence $w_\e$ which converges weakly to $w$ in $H(\Omega)$, we have to prove that
\begin{equation}\label{Ginf}
E(w)\leq \liminf_\e E_\e(w_\e) \,.
\end{equation}

For a given $w\in H(\Omega)$, let us first consider a regular sequence  $\widetilde{w}_\e \in C ^ \infty (\overline \Omega) \cap H (\Omega)$ weakly converging to $w$. 

As already done in the proof of Lemma \ref{propbounds} , by exploiting the $C^1$ regularity of the integral functionals $I_f$ and $I_g$ and the $L^2$-Mosco convergence of the corresponding sequences $\Delta_{\e, f}$ and $\Delta_{\e, g}$ defined according to \eqref{defDelta}, we infer
\begin{equation}\label{Ginf-regular}
E(w)\leq \liminf_\e E_\e(\widetilde{w}_\e)\,.
\end{equation}

Let us now consider a generic sequence ${w_\e}\in {H}(\Omega)$ weakly converging to $w$. By density, for every fixed $\e$, there exists $\widetilde{w}_{\e,k}\in C ^ \infty (\overline \Omega) \cap H (\Omega)$ such that  $\widetilde{w}_{\e,k}\to w_\e$ strongly in $H(\Omega)$ as $k\to +\infty$. 
With a diagonal argument, we may find a subsequence such that
\begin{equation}\label{diagonal}
\liminf_\e E_\e (w_\e) = \liminf_\e \lim_k E_\e (\widetilde{w}_{\e,k}) =\lim_\e  E_\e (\widetilde{w}_{\e,k_\e})\,,
\end{equation}
where the first equality follows by the strong continuity of $E_\e$.
We remark that the subsequence $\widetilde{w}_{\e,k_\e}$ is regular and weakly converging to $w$, hence it satisfies  (\ref{Ginf-regular}).  Such property, combined with the equality (\ref{diagonal}), concludes the proof of
(\ref{Ginf}).

\smallskip

Let us now prove (\ref{G_limsup}). For every fixed $w\in H(\Omega)$ we have to find a recovery sequence, namely a sequence $w_\e$ which converges weakly to $w$ in $H(\Omega)$ and satisfies
$
\limsup_{\e \to 0 }  E_\e(w_\e)\leq E(w) 
$.
If $\widetilde{w}$ is an element of the space $C ^ \infty (\overline \Omega) \cap H (\Omega)$, we are done simply by taking the constant sequence $w_\e\equiv \widetilde{w}$. Indeed, by following the same procedure adopted in the proof of Lemma \ref{propbounds}, we infer that
$
\lim_{\e\to 0} E_{\e}(\widetilde{w}) = E(\widetilde{w})
$.
If $w$ is a generic element of $H(\Omega)$, we approximate it by a sequence $\widetilde{w_k}\in C ^ \infty (\overline \Omega) \cap H (\Omega)$:  by the lower semicontinuity of the l.h.s.\ of  (\ref{G_limsup}) (usually called $\Gamma$-$\limsup E_\e$) and the continuity of $E$, we obtain  
$$
(\Gamma-\limsup E_\e )(w) \leq \liminf_k (\Gamma-\limsup E_\e )(\widetilde{w_k}) \leq \liminf_k \limsup_\e E_\e(\widetilde{w_k}) = E(w)\,.
$$

\qed

%
%
%


\medskip

\section{Regularity of the solution and necessary optimality condition}\label{secreg}

In this section we prove Propositions \ref{theo-reg} and \ref{theo-div}.
Both of them are based on the crucial result given below, which in turn exploits the $\Gamma$-convergence statement given in the previous section ({\it cf.}\ Lemma \ref{propgamma}).

\begin{proposition}\label{prop_w}
Under the standing assumptions, if $V$ has compact support contained into $\Omega$,  then the  function
$\overline \theta _V :=\langle V , \nabla \u \rangle $
solves the minimization problem
$\inf \limits_{w \in H ( \Omega)} E (\cdot, V)$, where $E (\cdot, V)$ is the functional  defined by \eqref{E}.
\end{proposition}

\proof
We first prove the following claim: for every fixed $\e >0$, the function
\begin{equation}\label{w_e} w_\e (x) := \frac{\u(x + \e V(x))-\u(x)}{\e}
\end{equation}
solves the minimization problem
$\inf \limits_{w \in H ( \Omega)} E_\e (\cdot, V)$, where $E_\e(\cdot, V)$ is the functional  defined by \eqref{E_e}.

\smallskip

Let $u_\e \in H(\Omega)$ be optimal for problem $J(\Omega_\e)$ in the formulation (\ref{Je2}) set on the fixed domain $\Omega$.
Via change of variables we obtain
$$
J(\Omega_\e) = - \int_{\Omega}[f_\e(\nabla u_\e) + g_\e(u_\e)]\,dx = - \int_{\Omega_\e} [f(\nabla (u_\e \circ \Psi_\e^{-1})) + g(u_\e\circ \Psi_\e^{-1})]\, dx\, .
$$
We infer that the function $u_\e \circ\Psi_\e^{-1} \in H(\Omega_\e)$ is optimal for problem $J(\Omega_\e)$ in its original formulation as an infimum over $H(\Omega_\e)$.
On the other hand, since $V$ has compact support contained into $\Omega$, we know that $J(\Omega_\e)=J(\Omega)$, and hence (since $\u$ is the  {\it unique} solution to $J (\Omega)$) we have $u_\e \circ \Psi_\e^{-1} = \u$, namely $u_\e (x)= \u \circ \Psi_\e(x) = \u (x+\e V(x))$.
It follows that  the function $w _\e$ defined in \eqref{w_e} minimizes $E_\e(\cdot, V)$
over $H (\Omega)$, and our claim on the optimality of $w_\e$ is proved.
\smallskip

In view of Lemma \ref{propgamma}, the minimizing sequence $\{w_\e\}$ is equibounded in $H(\Omega)$ and, up to subsequences, it weakly converges to
a minimizer of the $\Gamma$-limit,  namely to a minimizer of $E(\cdot, V)$. Clearly, the weak limit of $w_\e$ necessarily agrees with the pointwise limit $\langle V , \nabla \u \rangle $, which allows to conclude that
$\langle V,  \nabla \u \rangle$ minimizes $E (\cdot, V)$ over $H (\Omega)$.
\qed

\bigskip
{\it Proof of Proposition \ref{theo-reg}}.
In view of Proposition \ref{prop_w} we can assert that, for every vector field $V$ compactly supported into $\Omega$, the function $\langle V , \nabla \u \rangle $ belongs to $H ^ 1(\Omega)$.
By the arbitrariness of $V\in C^1_0(\Omega;\re^n)$, we infer that  $\u$ is an element of $W^{2,2}_{\rm loc}(\Omega)$.
\qed

\bigskip
{\it Proof of Proposition \ref{theo-div}}. 
We consider the minimization of the functional
$E (\cdot, V)$ in \eqref{E} over $H (\Omega)$, written as done in the proof of Proposition \ref{thm_J''}, see eq.\ \eqref{eq1}. By duality, we have the equality  \eqref{eq2}. Then, 
by   Lemma \ref{prop_duality},
 a function $w \in H (\Omega)$ and a pair $(\eta, \tau) \in L ^ 2 (\Omega; \re ^n) \times L ^ {q'} (\Omega)$ are optimal respectively for the two problems in \eqref{eq2} if and only if  
$(\eta, \tau) \in \partial \Psi (\nabla w, w)$ and $(\div \eta, \tau) \in \partial \Phi (w)$. 
Then, by computing explicitly the expressions of  $\partial \Phi$ and $\partial \Psi$, we get $\tau = \div \eta$, 
 with the additional condition $\eta \cdot n = 0 $ on $\Gamma _N$,
and
$$
\left\{
\begin{array}{ll}
\eta & = \nabla ^2 f(\nabla\u) (\nabla w  - DV ^ T \nabla \u)  - (DV-\div V\,I)\osigma \, ,
\\
\div\eta  &= g''(\u) w + g'(\u) \div V\, .
\end{array}
\right.
$$
Hence, for $V\in C^1_0(\Omega;\re^n)$, the distributional equality (\ref{div_eta}) follows by replacing
$w$ by  $\langle V , \nabla \u \rangle $  in the above optimality conditions (since we know by
Proposition \ref{prop_w} that $\overline \theta _V =\langle V , \nabla \u \rangle $ minimizes
$E (\cdot, V)$ over $H (\Omega)$).
For $V\in C^1(\re^n;\re^n)$,  the equality (\ref{div_eta})  continues to hold by the following simple argument:
in order to check \eqref{div_eta} for a 
given  test function $\varphi\in \mathcal D(\Omega)$, it is enough to 
use its validity when 
$V$ is replaced by the product $\mu V \in C^1_0(\Omega;\re^n)$, being $\mu$ a cut off function compactly supported into $\Omega$ which equals $1$ on
the support of $\varphi$.

\qed

\section{Representation formula for the second order shape derivative}\label{secrep}

In this section we give the

\bigskip
{\it Proof of Theorem \ref{bdry} (ii)}. 
We start from the following equality that we have already proved in Proposition \ref{thm_J''}:
\begin{equation}\label{J''=-inf}
J'' (\Omega, V) = - \inf _{w \in H (\Omega)}  E (w, V) \,.
\end{equation}
Throughout the proof we set for brevity $\w:= \langle V , \nabla \u \rangle $. 
By inserting into the definition \eqref{E} of $E (w, V)$ the identities
$$\begin{array}{ll}
& Q_f ( \nabla w- D V ^ T \nabla \u) = Q_f \big ( \nabla ( w - \w  ) \big ) +  \big\langle  \nabla ^ 2 f (\nabla \u) (\nabla ^ 2 \u) V , \nabla  (w -\w ) \big\rangle
+ Q _ f \big ( (\nabla ^ 2 \u ) V \big ) \\ \noalign{\medskip}
& Q_g (w) = Q_g ( w - \w  ) + g'' (\u) \,
\w \,   ( w - \w ) + Q _g (\w)\, ,
\end{array}$$
we obtain
$$\begin{array}{ll}& E(w, V)  \displaystyle =
 \int _{\Omega} 2 \big [   Q_f \big ( \nabla ( w - \w) \big ) +   Q_g ( w - \w )  \big ] \, dx \, + \\ \noalign{\medskip}
&\displaystyle+  \int_{\Omega}2  \big [ \langle \nabla ^ 2 f (\nabla \u) (\nabla ^ 2 \u) V - ( DV - \div V I ) \osigma , \nabla ( w - \w) \rangle + \big ( g'' (\u) \, \w +\div V g' (\u) \big ) ( w - \w)  \big ] \, dx \, + \\ \noalign{\medskip}
& \displaystyle
 + \int _\Omega \big [  \big ( f (\nabla \u) + g (\u)  \big )  (\div  V I -  DV ) : DV^T - 2 \langle (DV - \div V I) \osigma , (\nabla ^ 2 \u) V \rangle +2 \div V g' (\u) \w  \big ] \, dx \, + \\ \noalign{\medskip}
& \displaystyle
+ \int _\Omega \big [ 2 Q_f ((\nabla ^ 2 \u) V) + 2 Q_g (\w) \big ] \, dx \,.
\end{array}$$
Recalling the definition \eqref{defQ} of the quadratic form $\mathcal Q  ( \u, \cdot)  $ and 
the definition \eqref{defB} of the vector field $B (\u, V)$, by using also Proposition \ref{theo-div} 
we can rewrite the above equality as
$E(w, V) = {\rm (I)} + {\rm ( II)}$,
with
$${\rm (I)} = \mathcal Q ( \u,  w -\w ) \, dx   +  \int_{\Omega}2 \big [ \langle B (\u, V) ,  
\nabla ( w - \w) \rangle + 
\div B (\u, V) ( w - \w)  \big ] \, dx
$$ 
and
$$\begin{array}{ll} {\rm (II)} &=\displaystyle  
 \int _\Omega \big [  \big ( f (\nabla \u) + g (\u)  \big )  (\div  V I -  DV ) : DV^T   - 2 \langle (DV - \div V I) \osigma , (\nabla ^ 2 \u )V \rangle + 2 \div V g' (\u) \w \big ] \, dx \, + \\ \noalign{\medskip}
& \displaystyle
+ \int _\Omega \big [ 2 Q_f ((\nabla ^ 2 \u) V) + 2 Q_g (\w) \big ] \, dx \,.
\end{array}$$
Note that, in view of the assumption $\u\in W^{2,2}(\Omega)$, the vector field $B(\u,V)$ is in $L^2$ and has bounded divergence, thus it admits a normal trace in $H ^ {-1/2}(\partial \Omega)$,  and we are allowed to apply the integration by parts formula \eqref{teo-div}. Thus we get  
\begin{equation}\label{numero1} 
 {\rm (I)} =  \mathcal Q ( \u ,  w - \w ) \, dx + \int_{\partial \Omega} 2  (w-\w)   \, B (\u, V) \cdot n \, d \mathcal H ^ {n-1}\,, 
\end{equation}

We now want to rewrite (II) as a boundary integral. To this purpose we observe that, by  \eqref{defB} 
and \eqref{div_eta}, there hold
\begin{align}
-\langle (DV - \div V I ) \osigma - \nabla^2 f (\nabla^2 \u)V , (\nabla^2\u)V \rangle& =  \langle B(\u,V) , (\nabla^2\u)V\rangle\,, \label{id1}
\\ \noalign{\medskip}
\div V g'(\u)\, \w  + g''(\u) (\w)^2 &=  \div B(\u,V) \, \w \, .\label{id2}
\end{align}
Moreover, again thanks to the assumption $\u\in W^{2,2}(\Omega)$, the vector field
\begin{equation}\label{defX1} 
X_1:= (f(\nabla \u)+g(\u)) (\div V I - DV) V
\end{equation}
is bounded and has $L^2$ divergence. Thus we are allowed to apply the integration by parts formula \eqref{teo-div} and we get  
\begin{align} \int_\Omega (f(\nabla \u)+g(\u)) (\div V I - DV):DV^T \,dx & = \int_{\partial \Omega} (f(\nabla \u)+g(\u)) \langle(\div V I - DV)V,n\rangle  \, d \mathcal H ^ {n-1} + \notag
 \\ \noalign{\medskip}  & - \int_\Omega \langle (\div V I - DV) V,  (\nabla^2 \u) \osigma + g'(\u) \nabla \u  \rangle\,dx\, .\label{id3}
\end{align}

Using \eqref{id1}, \eqref{id2} and \eqref{id3}, we infer that
$$
\begin{array}{ll}{\rm (II)}  & = \displaystyle 
\int_{\partial \Omega} (X_1\cdot n)\, d \mathcal H ^ {n-1} +
\int_\Omega  [  \langle B(\u,V) , (\nabla^2\u)V\rangle + \div B(\u,V) \, \w ]\,dx \,+
\\ \noalign{\medskip}
& \ \
\displaystyle
- \int_\Omega  \langle (\nabla^2\u) DV \osigma - DV^T (\nabla^2 \u) \osigma - g'(\u) DV^T \nabla \u, V \rangle \, dx 
\,.\end{array}
$$
Exploiting the equality $(\nabla^2\u)V = \nabla \w - DV^T\nabla \u$, we get
$$
\begin{array}{ll}
{\rm (II)} & = \displaystyle \int_{\partial \Omega} (X_1 \cdot n) \, d \mathcal H ^ {n-1} + \int_\Omega  [  \langle B(\u,V) , \nabla \w \rangle + \div B(\u,V) \, \w ]\,dx\, +
\\ \noalign{\medskip}
& \ \
\displaystyle
- \int_\Omega \Big[
\langle  \nabla ^2 f(\nabla\u) (\nabla ^ 2 \u) V , DV^T \nabla \u\rangle    +\div V\langle \osigma , DV^T \nabla \u\rangle
+\langle (\nabla^2\u) DV \osigma, V \rangle \Big]\,dx \,+
 \\ \noalign{\medskip}
& \ \
\displaystyle
+ \int_\Omega \Big [\langle DV\osigma , DV^T \nabla \u\rangle + \langle  DV^T (\nabla^2 \u) \osigma, V \rangle
 + g'(\u)\langle DV^T \nabla \u, V \rangle
\Big ] \, dx \,.\end{array}
$$
Finally we remark that
$\langle B(\u,V) , \nabla \w \rangle + \div B(\u,V)\, \w  = \div \big ( \w  B(\u,V) \big)$,
and
\begin{align}
& - \langle  \nabla ^2 f(\nabla\u) (\nabla ^ 2 \u) V , DV^T \nabla \u\rangle    - \div V\langle \osigma , DV^T \nabla \u\rangle
-\langle (\nabla^2\u) DV \osigma, V \rangle + \langle DV\osigma , DV^T \nabla \u\rangle + \notag
 \\ \noalign{\medskip}
& \ \
\displaystyle
+ \langle  DV^T (\nabla^2 \u) \osigma, V \rangle
+ g'(\u)\langle DV^T \nabla \u, V \rangle = \div \Big(-\langle DV \osigma, \nabla \u\rangle V + \langle DV\, V, \nabla \u \rangle \osigma \Big)\,.\notag
\end{align}
Hence
\begin{equation}\label{numero2}
{\rm (II)} =
 \int_{\partial \Omega} (X_1 \cdot n) \, d \mathcal H ^ {n-1}  + \int_\Omega
 \div \Big ( \w \,  B(\u,V)  - \langle DV \osigma, \nabla \u\rangle V + \langle DV\, V, \nabla \u \rangle \osigma
  \Big)\,dx 
\,;
\end{equation}
moreover, thanks to the assumption $\u \in W ^ {2,2} (\Omega)$, the tensor field
\begin{equation}\label{defX2}
X_2:=  \w \, B (\u, V)     - \langle DV \osigma, \nabla \u\rangle V + \langle DV\, V, \nabla \u \rangle \osigma\,,
 \end{equation} 
is in $L ^2$ and has $L ^2$ divergence. Then, by adding \eqref{numero1} and \eqref{numero2} and by  applying the integration by parts formula \eqref{teo-div}, we obtain 
\begin{equation}\label{Efinal}
\begin{array}{ll}
E (w, V) = & \displaystyle \mathcal Q ( \u,  w - \w ) \, dx   +  \int_{\partial \Omega}
2(w-\w )\, B (\u, V) \cdot n \, d \mathcal H ^ {n-1}  +
\\ 
\noalign{\smallskip}
 & \displaystyle + 
  \int_{\partial \Omega} (X_1 \cdot n) \, d \mathcal H ^ {n-1}  +  \int_{\partial \Omega} (X_2 \cdot n) \, d \mathcal H ^ {n-1}  \,.\end{array}
\end{equation}

Now, the fact that the field $\mathcal C (\u, V)$ defined in (\ref{defC}) admits a normal trace in $H ^ {-1/2} (\partial \Omega)$ follows 
by observing that
$\mathcal C(\u,V) = - X_1 - X_2$, with $X_1$ and $X_2$ defined respectively 
in \eqref{defX1} and \eqref{defX2}. 
Finally 
the expression (\ref{J''1}) of the second order shape derivative follows by combining (\ref{Efinal}) 
with (\ref{J''=-inf}), and using the equality $\mathcal C(\u,V) = - X_1 - X_2$.

\qed

\section{Variants and perspectives}\label{secvar}

\subsection{The case $g$  linear}\label{subsec-linear}

We claim that Theorem \ref{bdry}, Proposition \ref{theo-reg}, and Proposition \ref{theo-div} remain true in the case when $g (v) = - \lambda v$,
for some $\lambda \in \re$. 
Note just that in the definition \eqref{defQ} of the quadratic form $\mathcal Q  (\u, \cdot)$, the term $Q_g$ vanishes. 

The proof of the lower bound inequality \eqref{thesis_lb} works unaltered,  whereas the proof of the upper bound inequality \eqref{thesis_ub} must be modified as follows. 
Note firstly that, since the Fenchel conjugate of $g$ is the Dirac delta function $g^*(\tau)=\delta_{-\lambda}$,
the dual energy $E ^ * (\eta, V)$ introduced in  \eqref{E*} is finite only if $\div \eta = - \lambda \div V$, and in this case it reads: 
\begin{equation}\label{E*bis}
 \begin{array}{ll}
 E ^ * (\eta, V)=&  \displaystyle 2 \int_\Omega [ f^*(\osigma)-   \langle \nabla f^*(\osigma),\osigma \rangle  + \lambda  \u ] \, a_2(DV) \,dx \,+ 
 \\
 &  \displaystyle  + 2  \int_\Omega \big [Q_{f^*}((DV - \div V \, I )\overline \sigma  + \eta)
+ \langle \nabla f^*(\osigma), DV\,\eta \rangle \big ]\,dx\, . 
\end{array}
\end{equation}
%

Next, 
we write the following inequality (which replaces \eqref{initial_ub1}): 
\begin{equation}\label{initial_ub2}
r_\e(V) \leq \frac{2}{\e^2} \left( \int_\Omega [f^*_\e (\sigma_\e) + g^*_\e (\div \sigma_\e) ]\,dx - J^*(\Omega,V) -\e J'(\Omega,V) \right) \, ,
\end{equation}
where $\sigma _\e$ are the perturbations of $\overline\sigma$ defined by
\begin{equation}\label{miss}
\sigma_\e:= \osigma + \sum_{k=1}^n \e^k \eta_k\, ,
\end{equation}
with for $k\in\{1,2,\dots ,n\}$
\begin{equation}\label{eta_k}
\eta_k\in   L^ \infty (\Omega;\re^n)  \, , \qquad 
\div \eta_k = - \lambda a_k(DV)\ \,\mathrm{in\ }\Omega\, ,\qquad \eta_k\cdot n =0 \ \mathrm{on\ }\Gamma_N\,.
\end{equation}

Notice that, for such a choice of  $\sigma _\e$, the r.h.s.\ of the inequality (\ref{initial_ub2}) is finite. Namely, 
exploiting the expansion (\ref{betae}) and the fact that  by \eqref{subdif} it holds $\div \osigma = -\lambda$, we obtain
$\div \sigma_\e = - \lambda \beta_\e$, so that 
\begin{equation}\label{g*=0}
g_\e^*(\div \sigma_\e)= \beta_\e g^*(\beta_\e^{-1} \div \sigma_\e) = \beta_\e g^*(-\lambda) = 0 \, . 
\end{equation}

Notice also that the set of fields satisfying (\ref{eta_k}) is non empty: indeed, every $a_k(DV)$ can be expressed as the divergence of a suitable regular vector field $\nu_k$ (see for instance \cite[Lemma 4.6.4]{M}); thus the conditions in (\ref{eta_k}) are verified by taking $\eta_k:= -\lambda \nu_k + \nabla w_k$, being $w_k$ the solution of
$$
\Delta w_k =0 \quad\text{in }\Omega\,,
\qquad \partial_n w_k = h_k \quad \text{on }\partial\Omega\,,
$$ 
with 
$
h_k= \lambda \nu_k \cdot n \ \text{on }\Gamma_N\,,\ \int_{\partial \Omega} h_k=0
$.
By combining (\ref{initial_ub2}) and (\ref{g*=0}), we may write
\begin{equation}\label{stima_ub2}
 r_\e(V) \leq  \Big (I_0^*(\e)  + I_1^* (\e) + I_2^*(\e)  \Big )\,,
\end{equation}
with
\begin{align*}
I_0^*(\e):= &  2 \int_\Omega f^*(\beta_\e^{-1} D\Psi_\e \sigma_\e)(a_2(DV) + \e m_\e) \,dx \, ,
\\
I_1^*(\e):= &  \frac{2}{\e} \int_\Omega[  f^*(\beta_\e^{-1}D\Psi_\e\sigma_\e)   - f^*(\osigma) ]\div V\,dx \, ,
\\
I_2^*(\e):= &  \frac{2}{\e^2} \int_\Omega [
f^*(\beta_\e^{-1}D\Psi_\e\sigma_\e) - f^*(\osigma)  - \e \langle \nabla f^*(\osigma), (DV- \div V I) \osigma \rangle - \e \lambda \div V \u ]\,dx  \, .
\end{align*}
Now we argue in a similar way as in the strictly convex case: we apply Lemma \ref{lemma_taylor} to the integral functional $I_{f^*}$ and, in order to deal with $I_2^*(\e)$, we exploit the equality
$$
\int_\Omega \lambda \div V \u \, dx =  - \int_\Omega \u \div \eta_1\,dx = \int_\Omega \langle \nabla f^*(\osigma) , \eta_1\rangle\,dx\, ;
$$
thus we obtain 
\begin{align*}
& \lim_{\e\to 0} I_0^*(\e) =  2 \int_\Omega f^*(\osigma )a_2(DV) \,dx\,,
\\
& \lim_{\e\to 0} I_1^*(\e) =  2 \int_\Omega \langle \nabla f^*(\osigma), (DV-\div V\,I)\osigma + \eta_1 \rangle \div V \,dx\,,
\\
& \lim_{\e\to 0} I_2^*(\e) =  2 \int_\Omega Q_{f^*}((DV-\div V\,I)\osigma + \eta_1 ) +
\\
& \quad \quad \quad \quad\quad\ \  \ \ + \langle \nabla f^*(\osigma),  ((\div V)^2 - \div V \,DV - a_2(DV))\osigma  + (DV-\div V\,I)\eta_1 + \eta_2      \rangle \,dx\,.
\end{align*}

In view of (\ref{stima_ub2}), by adding up the three terms above and exploiting the equality
$
\int_\Omega \langle  \nabla f^*(\osigma) , \eta_2 \rangle \,dx =  \int_\Omega \lambda \u a_2(DV)\,dx
$,
we infer that
$\limsup _{\e \to 0 } r _\e (V) \leq E ^*(\eta_1, V)$.

Now, since $\eta _1$ is an arbitrary vector field in $L ^ \infty (\Omega; \re ^n)$
with $\div \eta _1 = - \lambda$ and $\eta _1 \cdot n =  0$ on $\Gamma _N$, in order to conclude the proof  of the upper bound inequality
\eqref{thesis_ub} it is enough to show that every field $\eta \in X (\Omega; \re ^n)$ satisfying the two conditions $\div \eta = - \lambda$ and $\eta \cdot n =  0$ on $\Gamma _N$  can be approximated strongly in $L ^2$ by a sequence of fields $\eta ^h$ belonging to $L ^ \infty (\Omega; \re ^n)$ and satisfying the same two conditions. 
To that aim we begin by noticing that, since by assumption $\partial \Omega$ is piecewise $C^1$, there exists a sequence 
of smooth fields $\tilde \eta ^h$, with $\tilde \eta ^h \cdot n =  0$ on $\Gamma _N$, 
such that $\tilde \eta ^h \to \eta$ in $L ^ 2 (\Omega; \re ^n)$ and $\div \tilde \eta ^h \to  - \lambda$ in $L ^ 2 (\Omega)$. 
Then, we let $\phi ^h$ be the unique solution to the Dirichlet-Neumann problem
$$
- \Delta \phi ^ h = \rho ^h :=\lambda + \div \tilde \eta ^h  \quad  \text{in }\Omega\, , \qquad
\phi ^h = 0\quad \text{on }\Gamma _D\, , \qquad 
\partial_n\phi ^h = 0\quad  \text{on }\Gamma _N\,. 
$$ 
Finally we define $\eta ^h := \tilde \eta ^h + \nabla \phi ^h$, and we claim that such sequence $\eta ^h$ has the required properties. Namely, since $\rho ^h$ is a sequence of smooth functions converging to $0 $ in $L ^ 2 (\Omega)$, 
the sequence   $\eta ^h$ lies in $L ^ \infty (\Omega; \re ^n)$, and converges to $\eta$ in $L ^ 2(\Omega; \re ^n)$; moreover it holds $\div \eta^h = - \lambda$ and $\eta^h \cdot n =  0$ on $\Gamma _N$. 

Eventually, once now \eqref{thesis_ub} is established,  the remaining parts of the proofs given in Sections \ref{secexi}, \ref{secreg} and \ref{secrep} can be repeated as done for $g$ strictly convex. 

\subsection {The $p$-torsion problem}\label{sec-p}

Here we prove formula (\ref{formulap}) for the second order shape derivative of the $p$-torsional rigidity functional $J _p (\Omega)$ introduced in (\ref{def-p}), for $p > 2$.    
We recall that the first order derivative $J_p'(\Omega, V)$ exists and is given by 

$$ J_p'(\Omega,V) =  
\int_{\partial \Omega}  \frac{|\nabla \u|^p}{p'}\, V_n\, d\mathcal H^{n-1}\,,
$$
where $\u$ is the unique solution in $W ^ {1,p}_ 0(\Omega)$ to $- \Delta _p \u =\lambda$ (see for instance \cite{BFL, ChMa}).  It is well-known that  $\u$ is of class $C ^ {1, \alpha}(\Omega)$ (see \cite{DiBe, T}), and it is also in $C ^ 2 (\Omega \setminus S)$, where $S$ is the critical set $S:=\{x\in \Omega \, :\, \nabla \u = 0 \}$. Let us also recall that $S$  has vanishing Lebesgue measure (see \cite{Lo}), and it is compactly contained into $\Omega$ (thanks to Hopf boundary lemma).

Let us introduce the weight function $\rho:= |\nabla \u | ^ {p-2}$, which is continuous in $\Omega$ and  strictly positive outside $S$. 
Moreover, $\rho ^ {-1}$ is in $L ^ 1 (\Omega)$ (see \cite{Sci}). 

We denote by $W^{1, 2} _\rho(\Omega)$ the Hilbert space consisting of functions $v \in W ^ {1, 1}_{
\rm loc} (\Omega)$ such that
$$\| v \| _{ W ^ {1, 2} _\rho(\Omega)}:=
 \|v\| _{L ^ 2 (\Omega)} ^ 2 +  \|\nabla v\| _{L ^ 2 _ \rho(\Omega)} ^2< + \infty\, , $$
and by $H^{1, 2} _\rho(\Omega)$  the completion of $C ^ 1 (\overline \Omega)$ with respect to the above norm.

Let $\mathcal Q (\u, \cdot)$ denote the quadratic functional   defined on $W^{1,2} _\rho(\Omega)$ by  
\begin{equation}\label{defP} 
\mathcal Q (\u, v) := \int _\Omega \langle P (\u) \nabla v  , \nabla v \rangle \, dx \,, \  \hbox{ with } \  P(\u) :=  |\nabla\u|^{p-2}\left( I + (p-2) \frac{\nabla \u}{|\nabla\u|} \otimes \frac{\nabla \u}{|\nabla\u|}\right)\,. \
\end{equation}
Notice that, since the norm of the matrix $P (\u)$  is controlled by $\rho$, the functional $\mathcal Q (\u, \cdot)$ is continuous on the Hilbert space $W ^ {1, 2} _\rho (\Omega)$. 
By adapting the approach developed in Section \ref{secexi}, we obtain the following second order differentiability result. 

\begin{theorem} \label{prop-p}
Let $\partial \Omega$ be of class $C^2$. Let $p>2$ and assume that the weight $\rho:= |\nabla \u| ^ {p-2}$ is such that \begin{equation}\label{HW}
H ^ {1, 2} _\rho(\Omega) = W^{1, 2} _\rho(\Omega)\,.
\end{equation}

Then, the functional $J_p(\Omega)$ is twice differentiable at $\Omega$ in any direction $V$, and 
for $V$  normal to the boundary it holds
\begin{equation}\label{p-torsion}
 J''_p (\Omega, V)=- \frac{1}{p} \int_{\partial \Omega}  V_n ^ 2 \big ( p \lambda \partial_n \u  + |\partial_n \u| ^ p H _{\partial \Omega} \big ) \, d \mathcal H ^ {n-1}  
- \inf_{ v \in W^{1,2}_\rho (\Omega)\atop v= - V_n \partial _n \u \text { on } \partial \Omega }
\mathcal Q (\u, v)\,.
\end{equation}
\end{theorem}

\begin{remark}  {\rm Notice that, since $\rho$ is strictly positive in a neighborhood of $\partial \Omega$, 
the space $ W ^ {1, 2} _\rho (\Omega)$ is embedded into $H ^ 1 _{\rm loc} (\Omega \setminus S)$, 
and hence the trace operator is well defined from $ W ^ {1, 2} _\rho (\Omega)$ into $H ^ {1/2} (\partial \Omega)$. In particular, the infimum in \eqref{p-torsion} is well-defined since, by the $C^2$ regularity of $\u$ in $\Omega \setminus S$, the function $\langle V, \nabla \u \rangle$ belongs to $H ^ 1 _{\rm loc} (\Omega \setminus S)$. Moreover, the infimum can be equivalently taken in $H ^ 1 (\Omega)$ (as written in Example \ref{examplep}). 
}

\end{remark}

\begin{remark}  {\rm Let us add a few comments on assumption \eqref{HW} and related bibliographical references.  
As general texts on weighted Sobolev spaces, we refer to \cite{K, N}. Moreover, we address to the paper \cite{CSC} for a counterexample showing that the $L ^1$-summability property of $\rho ^ {-1}$ is not strong enough to ensure the validity of \eqref{HW}. 
A sufficient condition would be $\rho \in \mathcal A_2$, where $\mathcal A_2$ is the Muckenhoupt class of functions satisfying 
$\sup _B \Big ( \frac{1}{|B|} \int _B \rho \Big )  \Big ( \frac{1}{|B|} \int _B \rho ^ {-1} \Big )  < + \infty $,
where $B$ varies among balls in $\re^n$ (see for instance \cite{Kil}), but such kind of regularity result seems hard to prove.  
On the other hand, equality \eqref{HW} holds true if the set $S$ has vanishing capacity in $H ^ {1,2} _\rho (\Omega)$, namely if there exists a sequence  of smooth functions
$\alpha _\e:\Omega \to [0, 1]$  which are equal to $1$ in a $\e$-neighborhood of $S$ and satisfy $\lim _\e \int_{\Omega} |\nabla \alpha _\e |  ^2 \rho \, dx = 0$. Indeed in this case, given $w \in 
W^ {1,2}_{\rho} (\Omega)$ with $ w =  - \overline \theta _V$ on $\partial \Omega$,
it is possible to construct an approximating sequence $\{w_\e \}  \in H ^ {1, 2}_{\rho} (\Omega)$ with $ w _\e=  - \overline \theta _V$ on $\partial \Omega$: it is enough to take 
 $w _\e :=  (1 - \alpha_\e) w$ (which belong to $H ^ 1 (\Omega)$, and hence to $H ^ {1,2} _\rho (\Omega)$). 
In particular, our assumption \eqref{HW} turns out to be satisfied whenever $S$ is a singleton. 
Concerning the geometry of $S$ recall that, if $\Omega$ is convex, 
the function $\u$ is power concave \cite{Sa}, and hence $S$ agrees with the set where $\u$ assumes its maximum. 
If in addition $\Omega$ is {\it strictly} convex,  it is likely true that $S$ is reduced to a singleton
({\it cf.}\ \cite{BG, KL}), or at least that it is of vanishing capacity, so that  \eqref{HW} holds true. 
} 
\end{remark}

{\it Proof of Theorem \ref{prop-p}.}
Since shape derivatives only depend on the behavior of the deformation field on the boundary, 
with no loss of generality we may choose $V$ vanishing in a neighborhood of $S$ (recall that $S \subset \subset \Omega$). 
Then, in order to show the equality \eqref{p-torsion}, we proceed along the same scheme adopted in Section \ref{secexi}, namely we show that  the lower and upper limits as $\e \to 0$ of the differential quotients $r_\e(V)$ defined in \eqref{r_e} are bounded respectively from below and from above by  the same quantity, which is precisely the r.h.s.\ of formula \eqref{p-torsion}. 
For convenience, we divide the remaining of the proof in three steps. 

\smallskip
{\it Step 1 (lower bound)}:  Thanks to the assumption $p\geq 2$, the following lower bound can be achieved by arguing exactly as done in the proof the inequality \eqref{thesis_lb} in Lemma \ref{propbounds}:

$$
\liminf _{\e \to 0} r _\e (V) \geq m_\infty := -\inf_{w\in C^{\infty}_0(\Omega)} E(w,V)\,,
$$
where $E (w, V)$ is defined according to \eqref{E} by taking therein $f(z) = |z|^p/p$, and $g (v) = - \lambda v$.

We claim that the field $B(\u,V)$ introduced in (\ref{defB}) still satisfies Proposition \ref{theo-div}, namely
\begin{equation}\label{div-p}
\div B (\u,V) = -\lambda \div V\quad \hbox{in\ } \mathcal D'(\Omega)\,.
\end{equation}
The proof of the above equality cannot be repeated as done in Section \ref{secreg}, because  the minimizer $\u$ of the functional $J _p (\Omega)$ satisfies the regularity condition 
$\u \in W ^ {2,2} 
(\Omega)$ only for $p \in (1,3)$;   
nevertheless, for $p \geq 3$, one has that $\u \in W ^ {2, q}
(\Omega)$ for any $q < (p-1)/(p-2)$ (see \cite[Proposition 2.2]{DaSc}). Thus  
$\nabla \u$ is always in $W ^ {1,1}
(\Omega)$, and 
the distributional Hessian of $\u$ appearing in the definition $B(\u, V)$ is well-defined as a function in $L ^ 1 
(\Omega)$. Moreover, the function $\u$ satisfies the crucial regularity condition
$|\nabla \u|^{p-2}\nabla \u \in W^{1,2}(\Omega)$  (see \cite[Corollary 2.1]{Sc}),
which enables us to obtain the proof of \eqref{div-p} as follows. 
Recall that $ \osigma = \nabla f (\nabla \u) = |\nabla \u|^{p-2}\nabla \u$; then, since $\nabla f (\cdot)$ is a locally Lipschitz function and $\nabla \u \in L ^ \infty (\Omega) \cap W ^ {1,1} 
(\Omega)$,  
we can apply the  chain rule and we get
$
 D \osigma = \nabla^2 f (\nabla \u) \nabla ^2 \u$ a.e.\ in  $\Omega$.
Thus the field $B(\u,V)$ may be rewritten in terms of $\osigma$ as
$
B(\u,V) = D \osigma V - (DV - \div V I) \osigma
$
and in particular it belongs to $L^2
(\Omega;\re^n)$. It is then straightforward to obtain \eqref{div-p}. 
Namely,  for any test function $\varphi \in C^\infty_c(\Omega)$, it holds 
\begin{align*}
 \langle \div (D\osigma V ) , \varphi \rangle & = - \langle D \osigma V , \nabla \varphi\rangle 
 = - \langle\partial_k \osigma_i V_k , \partial_i \varphi \rangle =  \langle \osigma_i \partial_k V_k , \partial_i \varphi \rangle +  \langle \osigma_i V_k , \partial^2_{ik} \varphi \rangle
\\ 
& =   \langle  - \div (\osigma \div V ) + \div (V \div \osigma) + \div (DV \osigma), \varphi \rangle 
\\
& =    \langle  \div [ (DV - \div V I ) \osigma ], \varphi \rangle - \langle \lambda \div V , \varphi \rangle\,,
\end{align*}
where we have used the equality 
$- \div \osigma= \lambda$ (holding by \eqref{subdif} applied with $g(t)= - \lambda t$). 

Now, we are in a position to rewrite $m_\infty$ 
similarly as in Theorem \ref{bdry} (ii). Indeed, one can check that the proof of Theorem \ref{bdry} (ii) given in Section  \ref{secrep} continues to work, thanks in particular to the identity \eqref{div-p} and to the $W ^ {2,q}$ regularity of $\u$ for $q < (p-1)/(p-2)$).  Thus, defining
$\mathcal C_D(\u, V)$ as in (\ref{CD}) (with $f(z) = |z|^p/p$), letting
$\mathcal Q(\u, \cdot )$ be given by \eqref{defP}, and setting  $\overline \theta _V = \langle V, \nabla \u \rangle$,
we obtain
\begin{equation}\label{p-lb}
\begin{array}{ll}
\liminf _{\e \to 0} r _\e (V) \geq m_\infty & \displaystyle = \int_{\partial \Omega} \mathcal C_D (\u,V)\,d\mathcal H^{n-1} - \inf_{w\in C ^ \infty_0(\Omega)} \mathcal Q( \u, w- \overline \theta _V)
\\ \noalign{\medskip}
 & \displaystyle = \int_{\partial \Omega} \mathcal C_D (\u,V)\,d\mathcal H^{n-1} - \inf_{w\in H ^ {1, 2} _{\rho, 0} (\Omega)} \mathcal Q( \u, w- \overline \theta _V)\,;
\end{array}
\end{equation} 
here  $H ^ {1, 2} _{\rho, 0} (\Omega)$ denotes the space of functions in $H ^ {1, 2} _\rho(\Omega)$ vanishing at $\partial \Omega$, and the last equality follows from 
the continuity of $\mathcal Q (\u, \cdot)$ on $ W ^ {1, 2} _\rho$.

\medskip
{\it Step 2 (upper bound)}: 
Due to the growth of order $p'\in (1,2]$ of the Fenchel conjugate $f^*(z^*) = |z^*| ^ {p'} /{p'}$, in order bound from above $r_\e (V)$, we need to set up an approximation argument.  Roughly speaking, we cover almost all $\Omega$ with a suitable family of subdomains of $\Omega$ in which the Hessian of $f^*$ is bounded from above, and then we pass to the limit. 
More precisely, as $\nabla \u$ is a continuous function which is strictly positive outside $S$,  we can construct  increasing sequence  of open sets $\Omega _h\uparrow (\Omega\setminus S) $ 
such that, for every $h \in \mathbb N$, $|\nabla \u| > 1/h$ in $\Omega _h$, and the boundary $\Gamma_h:=\partial \Omega_h \setminus \partial \Omega$ is smooth and has a positive distance from $\partial\Omega$. 
For every fixed $h$, we follow the same procedure adopted for the proof of the inequality \eqref{thesis_ub} in Lemma \ref{propbounds}, in the variant when $g$ is linear described in \S \ref{subsec-linear}. 
We choose the fields $\eta _k$ appearing in \eqref{miss} as done  in \eqref{eta_k}, with the additional condition that $\eta_k = 0$ in $\Omega \setminus \Omega _k$ 
(note that this is possible
thanks to the assumption made on the support of the deformation field $V$). 
We infer that, for every fixed $h$, 
\begin{equation}\label{e1}
\limsup_{\e} r_\e(V) \leq m_h:= \inf_{\eta\in X(\Omega;\svre^n)	\,,	\,	\, \div \eta = -\lambda \div V\, \mathrm{in\,} \Omega
\atop \qquad \eta = 0 \,\mathrm{in\,}\Omega \setminus \Omega_h 
}  E^*(\eta,V)\,,
\end{equation}
where $E ^* (\eta, V)$ is defined according to \eqref{E*bis} (with $f^*(z^*) = |z^*|^{p'}/{p'}$).


By applying Lemma \ref{prop_duality} and exploiting the regularity of $\u$, we may rewrite $m_h$ in primal form as 
\begin{equation}\label{e2}
m_h= \int_{\partial \Omega} \mathcal C_D (\u,V)\, d\mathcal H^{n-1} - \inf_{w\in V_h}\mathcal Q_h (\u, w-\overline \theta _V)\,,
\end{equation}
where $V_h$ denotes the subspace of $H ^1 (\Omega _h)$ of functions having zero trace on $\partial \Omega$, and 
$$
\mathcal Q_h(\overline u, v):= 2 \int_{\Omega_h} \langle P (\u) \nabla v , \nabla v \rangle  \,dx\,. 
$$
Now, in order to bound from above the upper limit of $m _h$ as $h \to + \infty$, we are going to bound from below  the infimum appearing at the r.h.s.\ of \eqref{e2}. We start by noticing that such infimum can be equivalently taken in the space 
$W ^ {1, 2} _{\rho, 0} (\Omega)$ of functions in $W ^ {1, 2} _\rho(\Omega)$ vanishing at $\partial \Omega$. 
Let  $\{ w _h \} \subset W ^ {1, 2} _{\rho, 0} (\Omega)$ be a sequence of minimizers for $\mathcal Q_h(\u, \cdot - \w)$. For every fixed $k \in \mathbb N$, such sequence is bounded in $H ^ 1 (\Omega _k)$, so that it admits a subsequence  weakly converging in $H ^ 1 (\Omega_k)$ to some element $w ^ {(k)}$. 
By a diagonalization argument, we can choose the same subsequence for every $k$, and hence the restriction of $w^{(k)}$ to $\Omega_l$ for $l<k$ agrees with $w^{(l)}$. 
Hence we obtain a (not relabeled) subsequence of $w _h$ and an element  $\hat w \in W ^ {1, 2} _{\rho, 0} (\Omega)$ such that
$w_h$ converges weakly to $\hat w$ in $H^1(\Omega_k)$ for every $k$. 
Then, by Fatou's Lemma and monotone convergence, we get
\begin{equation}\label{e3} 
\liminf _h \mathcal Q_h(\u,  w_h -\overline \theta _V) \ge \mathcal Q(\u, \hat  w -\overline \theta _V) \ge \inf _{w \in W ^ {1, 2} _{\rho, 0} (\Omega)} \mathcal Q(\u,  w -\overline \theta _V) \ .
\end{equation}
By combining \eqref{e1}, \eqref{e2}, and \eqref{e3}, we obtain
\begin{equation}\label{p-ub}
\limsup_{\e} r_\e(V) \leq m _\infty ^* := 
\int_{\partial \Omega} \mathcal C_D (\u,V)\,d\mathcal H^{n-1} -\inf _{w \in W ^ {1, 2} _{\rho, 0} (\Omega)} \mathcal Q(\u,  w -\overline \theta _V)\,. 
\end{equation}

\smallskip
{\it Step 3 (conclusion)}: We finally need to show that the lower and upper bounds 
$m _\infty$  and  $m_{\infty }^*$  in 
\eqref{p-lb} and \eqref{p-ub} agree, and that they are equal to the expression at the r.h.s. of \eqref{p-torsion}. 
We firstly observe that, by arguing as in Remark \ref{repre},  the integral over $\partial \Omega$ of $\mathcal C_D (\u,V)$ appearing in the expression \eqref{p-lb} of $m _\infty$ can be rewritten as done in \eqref{p-torsion} .
Then it only remains to prove that  
no Lavrenteev phenomenon occurs for the infimum problems
appearing in \eqref{p-lb} and \eqref{p-ub}. In other words, 
after a translation, we are reduced to show that 
\begin{equation}\label{Lav2} \inf_{w\in W^ {1,2}_{\rho} (\Omega)}  \left\{ \mathcal Q (\u, \cdot) \, : \, \hbox{$ w =  - \overline \theta _V$ on $\partial \Omega$ } \right\}
\ =\    \inf_{w\in  H ^ {1, 2}_{\rho} (\Omega)}   \left\{ \mathcal Q (\u, \cdot) \, : \, \hbox{$ w =  - \overline \theta _V$ on $\partial \Omega$ } \right\}\ .
\end{equation}

The validity of \eqref{Lav2} is an immediate consequence of our assumption \eqref{HW}
and of the continuity of $\mathcal Q (\u, \cdot)$ in $W ^ {1, 2} _\rho (\Omega)$.

\qed

\subsection {Perspectives} 

A natural perspective is the possible extension of our results  to the case of convex {\it non-smooth} integrands $f$ and $g$. 
The main difficulty is that in this case one can no longer exploit the second order differentiability property for integral functionals stated in Lemma \ref{lemma_taylor}. 
However we expect that in some cases the second order shape derivative still exists, and admits a representation formula similar to \eqref{J''1}, with  the quadratic form $\mathcal Q (\u, \cdot)  $ (formerly associated with the second order differentials of $f$ and $g$) replaced by a suitable convex, positively $2$-homogeneous, non-quadratic function.  
This kind of generalization would be useful in applications: for instance, the shape functional obtained by taking 
the function $
f(z)$ equal to $\frac{|z| ^ 2}{2} + \frac{1}{2}$  if  $|z| \geq 1$ and 
$|z|$ if  $|z| < 1$,  and the function $g (v )$ equal to $- \lambda v$   (for some $\lambda \in \re$), 
is related to the optimization of thin rods in torsion regime (see \cite{ABFL, BFLS}), and studying its second order shape derivative might be helpful in order to investigate the occurrence of homogenization regions in an optimal design. (For the computation of the first order shape derivative in this case, see Example 3.9 (ii) in \cite{BFL}.) 

\medskip

Another perspective is trying to understand whether the results obtained in \S \ref{sec-p} may entail some useful information on the second order shape derivative of  the shape functional $J _\infty (\Omega)$ obtained as the limit  as $p \to + \infty$ of the $p$-torsional rigidity functionals $J _p (\Omega)$ defined in (\ref{def-p}). Actually, up to constant multiple, it is well known that here holds
$J _\infty (\Omega) = \int _\Omega d_{\partial \Omega} (x) \, dx$ and $\lim _p u _p (x) = d _{\partial \Omega} (x) $,
where $d _{\partial \Omega}$ is the distance function from $\partial \Omega$, and $u _p$ is the unique solution to $J _p (\Omega)$ (see  \cite{BDM, BBD, PP}). 
%

\section{Appendix}\label{secapp}

\begin{lemma}\label{lemma_taylor}
Let $\phi:\re^n \to \re$ be a strongly convex function of class $C^2$ satisfying  growth conditions of order $2$ from above and below, and let $I _\phi$ be the integral functional defined on $L^2(\Omega;\re^n)$ by
\begin{equation}\label{Iphi}
 I_\phi(z):=\int_\Omega \phi(z(x))\,dx\,. \end{equation}
Let $z_0$ be a fixed vector field in $L^\infty(\Omega;\re^n)$ and, 
for $h \in L^2(\Omega;\re^n)$, set
$q_\phi(h):=\frac{1}{2}\, \int_\Omega \langle \nabla^2 \phi(z_0) h, h \rangle \, dx$, and 
\begin{equation}\label{defDelta} 
 \Delta_{\e,\phi}(h):=\frac{I_\phi(z_0 + \e h) - I_\phi(z_0) - \e \langle \nabla \phi(z_0) , h\rangle_{L^2}}{\e^2}\,.
\end{equation}
Then: 
\begin{itemize}
\item[(i)] $I_\phi$ is of class $C^1$ on the Hilbert space $L^2(\Omega;\re^n)$, and
$\nabla I _\phi (z) = \nabla \phi( z)$; 
\item[(ii)] $\Delta_{\e,\phi}$ Mosco converges in $L^2(\Omega;\re^n)$ to $q_{\phi}$ as $\e\to 0$;
\item[(iii)]  $I_\phi$ is second order differentiable at $z_0$ with respect to $L^\infty$ variations; namely, for every  $h\in L^\infty(\Omega;\re^n)$, it holds
$
\Delta_{\e, \phi}(h) \to q_{\phi}(h)$  as $\e \to 0$.
\end{itemize}
\end{lemma}

\proof
For statements  (i) and (ii), see  \cite[Proposition 4.1 and Corollary 3.2]{Noll1}. 
Let us now prove (iii).
For a fixed $h\in L^\infty(\Omega;\re^n)$, the Mosco convergence property stated in (ii) implies that
$$
q_{\phi}(h)\leq  \liminf\limits_{\e \to 0} \Delta_{\e,\phi} (h)\,.
$$ 
In order to prove that 
$\limsup_\e  \Delta_{\e,\phi} (h) \leq q_{\phi}(h)$,
we exploit the second property given by the Mosco convergence: there exists a sequence $h_\e\in L^2(\Omega;\re^n)$ such that $h_\e \stackrel{L^2}{\to} h$ and
$
q_{\phi}(h_\e)= \lim_{\e\to 0} \Delta_{\e,\phi} (h_\e)$. 
We remark that the function $ \Delta_{\e,\phi} $ is convex, therefore,
$
 \Delta_{\e,\phi} (h) - \Delta_{\e,\phi} (h_\e) \leq \langle \nabla \Delta_{\e,\phi} (h) , h-h_\e\rangle_{L^2}$.
If we prove that, for $\e>0$ small enough,
\begin{equation}\label{claim}
\left\| \nabla \Delta_{\e,\phi} (h) \right\|_{L^2} \leq C\, ,
\end{equation}
we are done, since we infer
$$
\limsup_\e  \Delta_{\e,\phi} (h) - q_{\phi}(h) = \limsup_\e \Big( \Delta_{\e,\phi} (h) -  \Delta_{\e,\phi} (h_\e) \Big)\leq \lim_{\e\to 0} C  \|h_\e - h\|_{L^2}=0\,.
$$
We conclude by proving the claim (\ref{claim}): we remark that
$
\nabla \Delta_\e (h) = [{\nabla \phi(z_0 + \e h) - \nabla \phi(z_0)}]/{\e}
$,
and hence \eqref{claim} follows from the local boundedness of the Hessian matrix $\nabla^2\phi$ and recalling that by assumption there exists $R>0$ such that
$\|z_0 +  t h\|_{L^\infty}\leq R$ for $t$ small enough.
\qed

%
%
%
%
%
%

\begin{lemma}\label{lemma_duality}
Let $X$ be a normed vector space and let $X^*$ be its topological dual. We denote by
$\langle \cdot , \cdot \rangle$ the duality product. Let $h:X\to \re\cup\{+\infty\}$ be a proper function. Then, for every $a\in X$ and $b\in X^*$, there holds
$$
\left(h(\cdot - a)\right)^* (z^*)= h^*(z^*) + \langle a, z^*\rangle\,, \qquad 
\left(h(\cdot) + \langle b,\cdot \rangle \right)^* (z^*)= h^*(z^*-b)\, .$$
\end{lemma}
\proof The statement follows straightforward by using the definition of Fenchel conjugate. \qed

\begin{lemma}\label{lemma_quadratic}
Let $A$ be a positive definite $n\times n$ matrix and let $Q_A$ denote the associated quadratic form, defined as $Q_A(x):=\frac{1}{2} \langle Ax , x\rangle$. Then
$
\left(Q_A\right)^*=Q_{A^{-1}}
$.
\end{lemma}
\proof See \cite[Chapter III]{R}.\qed

\begin{lemma}\label{lemma_hessian}
Let $X$ be a Banach space and let $h:X\to \re$ be a strongly convex function of class $C^2$.  Then $\mathrm{dom} (h^*)$ has nonempty interior, $h^*$ is $C^2$ on the interior of
$\mathrm{dom} (h^*)$,
and $
x^*=\nabla h(x)$ with $x \in X$ implies $
x=\nabla h^*(x^*)$ and 
$
\nabla^2 h^*(x^*) = \left(\nabla^2 h(x)\right)^{-1}$. 
\end{lemma}
\proof See \cite[Proposition 10 in Section 2 of Chapter II]{Ek}.\qed

\begin{lemma}\label{prop_duality}
Let $Y, Z$ be Banach spaces. Let $A:Y \to Z$ be a linear operator
with dense domain $D(A)$. Let $\Phi: Y \to \re \cup \{+ \infty\}$ be
convex, and $\Psi: Z \to \re \cup \{ + \infty\}$ be convex lower
semicontinuous. Assume there exists $u _0 \in D (A)$ such that $\Phi(u_0) < + \infty$ and $\Psi$ is continuous at $A\,u_0$.  Let $Z^*$
denote the dual space of $Z$, $A^*$ the adjoint operator of $A$, and
$\Phi ^*$, $\Psi ^*$ the Fenchel conjugates of $\Phi$, $\Psi$. Then
\begin{equation}\label{dual}
- \inf _{u \in Y} \Big \{ \Psi (A\,u) + \Phi (u)  \Big \} = \inf
_{\sigma \in Z ^*} \Big \{ \Psi ^* (\sigma) + \Phi ^* (- A ^*\,
\sigma) \Big \}\, ,
\end{equation}
and the infimum at the right hand side is achieved.

Furthermore,  $\ov{u}$ and $\ov{\sigma}$ are optimal for the l.h.s.\ and the r.h.s.\ of $(\ref{dual})$ respectively, if and
only if  there holds $\ov{\sigma}\in \partial
\Psi(A\ov{u})$ and $-A^\ast \ov{\sigma} \in \partial \Phi (\ov{u})$.
\end{lemma}

\proof See \cite[Proposition 14]{Bo}. \qed

\bigskip

\end{document}